\documentclass{article}
\usepackage{amsmath,amssymb,amsthm,eucal,eufrak}
\usepackage{graphicx}

\newtheorem{theorem}{Theorem}
\newtheorem{lemma}{Lemma}
\newtheorem{proposition}{Proposition}
\newtheorem{definition}{Definition}
\newtheorem{corollary}{Corollary}
\newtheorem{example}{Example}
\newtheorem{remark}{Remark}

\numberwithin{theorem}{section}
\numberwithin{definition}{section}
\numberwithin{lemma}{section}
\numberwithin{corollary}{section}
\numberwithin{equation}{section}
\numberwithin{proposition}{section}
\numberwithin{example}{section}
\numberwithin{remark}{section}
\numberwithin{figure}{section}

\def\vtl{\vskip 1mm}

\def\es{\varnothing}
\def\SB{\subseteq}

\def\aa{\alpha}

\def\hh{\theta}

\def\b0{\boldsymbol 0}

\def\Na{\mathbb N}

\def\Zee{\mathbb Z}

\def\imp{\Rightarrow}

\def\eq{\Leftrightarrow}

\def\AAA{{\cal A}}
\def\FFF{{\cal F}}

\def\GGG{{\cal G}}
\def\HHH{{\cal H}}
\def\III{{\cal I}}
\def\JJJ{{\cal J}}

\def\ZZZ{{\cal Z}}

\def\PPP{{\cal P}}

\def\begeq{\begin{equation}}
\def\edeq{\end{equation}}

\def\roster{\begin{enumerate}}
\def\endroster{\end{enumerate}}

\makeindex

\begin{document}

\title{Partial cubes: structures, characterizations, and constructions}
	\author{Sergei~Ovchinnikov\\ Mathematics Department\\San Francisco State University\\San Francisco, CA 94132\\sergei@sfsu.edu}
\date{May 8, 2006}
\maketitle

\begin{abstract}\noindent
Partial cubes are isometric subgraphs of hypercubes. Structures on a graph defined by means of semicubes, and Djokovi\'{c}'s and Winkler's relations play an important role in the theory of partial cubes. These structures are employed in the paper to characterize bipartite graphs and partial cubes of arbitrary dimension. New characterizations are established and new proofs of some known results are given.

The operations of Cartesian product and pasting, and expansion and contraction processes are utilized in the paper to construct new partial cubes from old ones. In particular, the isometric and lattice dimensions of finite partial cubes obtained by means of these operations are calculated.
\end{abstract}

\vtl\noindent

\section{Introduction} \label{S:introduction}

A hypercube $\HHH(X)$ on a set $X$ is a graph which vertices are the finite subsets of $X$; two vertices are joined by an edge if they differ by a singleton. A partial cube is a graph that can be isometrically embedded into a hypercube.

\vtl
There are three general graph-theoretical structures that play a prominent role in the theory of partial cubes; namely, semicubes, Djokovi\'{c}'s relation $\hh$, and Winkler's relation $\Theta$. We use these structures, in particular, to characterize bipartite graphs and partial cubes. The characterization problem for partial cubes was considered as an important one and many characterizations are known. We list contributions in the chronological order: Djokovi\'{c}~\cite{dD73} (1973), Avis~\cite{dA81} (1981), Winkler~\cite{pW84} (1984), Roth and Winkler~\cite{rR86} (1986), Chepoi~\cite{vC88,vC94} (1988 and 1994). In the paper, we present new proofs for the results of Djokovi\'{c}~\cite{dD73}, Winkler~\cite{pW84}, and Chepoi~\cite{vC88}, and obtain two more characterizations of partial cubes.

\vtl
The paper is also concerned with some ways of constructing new partial cubes from old ones. Properties of subcubes, the Cartesian product of partial cubes, and expansion and contraction of a partial cube are investigated. We introduce a construction based on pasting two graphs together and show how new partial cubes can be obtained from old ones by pasting them together.

\vtl
The paper is organized as follows.

\vtl
Hypercubes and partial cubes are introduced in Section~\ref{S:hypercubes} together with two basic examples of infinite partial cubes. Vertex sets of partial cubes are described in terms of well graded families of finite sets.

\vtl
In Section~\ref{S:characterizations} we introduce the concepts of a semicube, Djokovi\'{c}'s $\hh$ and Winkler's $\Theta$ relations, and establish some of their properties. Bipartite graphs and partial cubes are characterized by means of these structures. One more characterization of partial cubes is obtained in Section~\ref{S:fundamental sets}, where so-called fundamental sets in a graph are introduced.

\vtl
The rest of the paper is devoted to constructions: subcubes and the Cartesian product (Section~\ref{S:subcubes Cartesian}), pasting (Section~\ref{S:pasting}), and expansions and contractions (Section~\ref{S:expansions}). We show that these constructions produce new partial cubes from old ones. Isometric and lattice dimensions of new partial cubes are calculated. These dimensions are introduced in Section~\ref{S:dimensions}.

\vtl
Few words about conventions used in the paper are in order. The {\sl sum {\rm(}disjoint union{\rm)}} $A+B$ of two sets $A$ and $B$ is the union
$$
(\{1\}\times A) \cup (\{2\}\times B).
$$
All graphs in the paper are simple undirected graphs. In the notation $G=(V,E)$, the symbol $V$ stands for the set of vertices of the graph $G$ and $E$ stands for its set of edges. By abuse of language, we often write $ab$ for an edge in a graph; if this is the case, $ab$ is an unordered pair of distinct vertices. We denote $\langle U\rangle$ the graph induced by the set of vertices $U\SB V$. If $G$ is a connected graph, then $d_G(a,b)$ stands for the distance between two vertices $a$ and $b$ of the graph $G$. Wherever it is clear from the context which graph is under consideration, we drop the subscript $G$ in $d_G(a,b)$. A subgraph $H\SB G$ is an {\sl isometric subgraph} if $d_H(a,b)=d_G(a,b)$ for all vertices $a$ and $b$ of $H$; it is {\sl convex} if any shortest path in $G$ between vertices of $H$ belongs to $H$.

\section{Hypercubes and partial cubes} \label{S:hypercubes}

Let $X$ be a set. We denote $\PPP_f(X)$ the set of all finite subsets of $X$.

\begin{definition} \label{D:hypercube}
{\rm A graph $\HHH(X)$ has the set $\PPP_f(X)$ as the set of its vertices; a pair of vertices $PQ$ is an edge of $\HHH(X)$ if the symmetric difference $P\Delta Q$ is a singleton. The graph $\HHH(X)$ is called the {\sl hypercube on $X$}~\cite{dD73}. If $X$ is a finite set of cardinality $n$, then the graph $\HHH(X)$ is the {\sl $n$-cube} $Q_n$. The {\sl dimension} of the hypercube $\HHH(X)$ is the cardinality of the set $X$.
}
\end{definition}

The shortest path distance $d(P,Q)$ on the hypercube $\HHH(X)$ is the {\sl Hamming distance} between sets $P$ and $Q$:
\begeq \label{Hamming distance}
d(P,Q)=|P\Delta Q|\quad\text{for $P,Q\in\PPP_f$.}
\edeq
The set $\PPP_f(X)$ is a metric space with the metric $d$.

\begin{definition} \label{D:partial cube}
{\rm
A graph $G$ is a {\sl partial cube} if it can be isometrically embedded into a hypercube $\HHH(X)$ for some set $X$. We often identify $G$ with its isometric image in the hypercube $\HHH(X)$, and say that $G$ is a {\sl partial cube on the set} $X$.
}
\end{definition}

{\begin{figure}[h!]
\centerline{\includegraphics{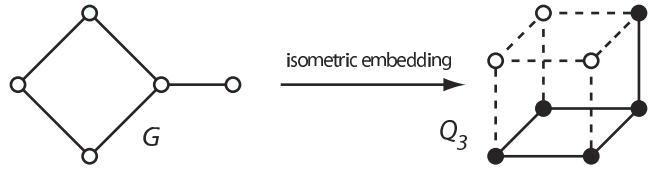}}
\caption{A graph and its isometric embedding into $Q_3$.} \label{example0} 
\end{figure}
}

An example of a partial cube and its isometric embedding into the cube $Q_3$ is shown in Figure~\ref{example0}.

\vtl
Clearly, a family $\FFF$ of finite subsets of $X$ induces a partial cube on $X$ if and only if for any two distinct subsets $P,Q\in\FFF$ there is a sequence
$$
R_0=P,R_1,\ldots,R_n=Q
$$
of sets in $\FFF$ such that
\begeq \label{wg-family}
d(R_i,R_{i+1})=1\quad\text{for all $0\leq i<n$,}\quad\text{and}\quad d(P,Q)=n.
\edeq

The families of sets satisfying condition~(\ref{wg-family}) are known as well graded families of sets~\cite{jDjF97}. Note that a sequence $(R_i)$ satisfying~(\ref{wg-family}) is a shortest path from $P$ to $Q$ in $\HHH(X)$ (and in the subgraph induced by $\FFF$).

\begin{definition} \label{D:wg-family}
{\rm
A family $\FFF$ of arbitrary subsets of $X$ is a {\sl wg-family} ({\sl well graded family of sets}) if, for any two distinct subsets $P,Q\in\FFF$, the set $P\Delta Q$ is finite and there is a sequence
$$
R_0=P,R_1,\ldots,R_n=Q
$$
of sets in $\FFF$ such that $|R_i\Delta R_{i+1}|=1$ for all $0\leq i<n$ and $|P\Delta Q|=n$.
}
\end{definition}

\begin{example} \label{nonisometric path}
{\rm The induced graph can be a partial cube on a different set if the family $\FFF$ is not well graded. Consider, for instance, the family
$$
\FFF=\{\es,\{a\},\{a,b\},\{a,b,c\},\{b,c\}\}
$$
of subsets of $X=\{a,b,c\}$. The graph induced by this family is a path of length $4$ in the cube $Q_3$ (cf. Figure~\ref{cube3}). Clearly, $\FFF$ is not well graded. On the other hand, as it can be easily seen, any path is a partial cube.
}
\end{example}

{\begin{figure}[h!]
\centerline{\includegraphics{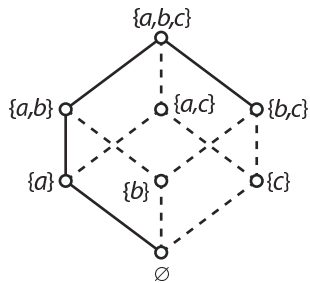}}
\caption{A nonisometric path in the cube $Q_3$.} \label{cube3} 
\end{figure}
}

Any family $\FFF$ of subsets of $X$ defines a graph $G_\FFF=(\FFF,E_\FFF)$, where 
$$
E_\FFF=\{\{P,Q\}\SB\FFF:|P\Delta Q|=1\}.
$$

\begin{theorem} \label{wg=finite wg}
The graph $G_\FFF$ defined by a family $\FFF$ of subsets of a set $X$ is isomorphic to a partial cube on $X$ if and only if the family $\FFF$ is well graded.
\end{theorem}

\begin{proof}
We need to prove sufficiency only. Let $S$ be a fixed set in $\FFF$. We define a mapping $f:\FFF\rightarrow\PPP_f(X)$ by $f(R)=R\Delta S$ for $R\in\FFF$. Then
$$
d(f(R),f(T))=|(R\Delta S)\Delta(T\Delta S)|=|R\Delta T|.
$$
Thus $f$ is an isometric embedding of $\FFF$ into $\PPP_f(X)$. Let $(R_i)$ be a sequence of sets in $\FFF$ such that $R_0=P$, $R_n=Q$, $|P\Delta Q|=n$, and $|R_i\Delta R_{i+1}|=1$ for all $0\leq i<n$. Then the sequence $(f(R_i))$ satisfies conditions~(\ref{wg-family}). The result follows.
\end{proof}

\vtl
A set $R\in\PPP_f(X)$ is said to be {\sl lattice between} sets $P,Q\in\PPP_f(X)$ if 
$$
P\cap Q\SB R\SB P\cup Q.
$$
It is {\sl metrically between} $P$ and $Q$ if
$$
d(P,R)+d(R,Q)=d(P,Q).
$$
The following theorem is a well-known result about these two betweenness relations on $\PPP_f(X)$ (see, for instance,~\cite{lB53}).

\begin{theorem} \label{lattice=metric betweenness}
Lattice and metric betweenness relations coincide on $\PPP_f(X)$.
\end{theorem}

Let $\FFF$ be a family of finite subsets of $X$. The set of all $R\in\FFF$ that are between $P,Q\in\FFF$ is the {\sl interval $\III(P,Q)$ between $P$ and $Q$ in $\FFF$}. Thus,
$$
\III(P,Q)=\FFF\cap[P\cap Q,P\cup Q],
$$
where $[P\cap Q,P\cup Q]$ is the usual interval in the lattice $\PPP_f$.

\vtl
Two distinct sets $P,Q\in\FFF$ are {\sl adjacent in $\FFF$} if $\JJJ(P,Q)=\{P,Q\}$. If sets $P$ and $Q$ form an edge in the graph induced by $\FFF$, then $P$ and $Q$ are adjacent in $\FFF$, but, generally speaking, not vice versa. For instance, in Example~\ref{nonisometric path}, the vertices $\es$ and $\{b,c\}$ are adjacent in $\FFF$ but do not define an edge in the induced graph (cf. Figure~\ref{cube3}).

\vtl
The following theorem is a `local' characterization of wg-families of sets.

\begin{theorem} \label{local wg}
A family $\FFF\SB\PPP_f(X)$ is well graded if and only if $d(P,Q)=1$ for any two sets $P$ and $Q$ that are adjacent in $\FFF$.
\end{theorem}

\begin{proof}
(Necessity.) Let $\FFF$ be a wg-family of sets. Suppose that $P$ and $Q$ are adjacent in $\FFF$. There is a sequence $R_0=P,R_1,\ldots,R_n=Q$ that satisfies conditions~(\ref{wg-family}). Since the sequence $(R_i)$ is a shortest path in $\FFF$, we have
$$
d(P,P_i)+d(P_i,Q)=d(P,Q)\quad\text{for all $0\leq i\leq n$.}
$$
Thus, $P_i\in\III(P,Q)=\{P,Q\}$. It follows that $d(P,Q)=n=1$.

\vtl
(Sufficiency.) Let $P$ and $Q$ be two distinct sets in $\FFF$. We prove by induction on $n=d(P,Q)$ that there is a sequence $(R_i)\in\FFF$ satisfying conditions~(\ref{wg-family}).

The statement is trivial for $n=1$. Suppose that $n>1$ and that the statement is true for all $k<n$. Let $P$ and $Q$ be two sets in $\FFF$ such that $d(P,Q)=n$. Since $d(P,Q)>1$, the sets $P$ and $Q$ are not adjacent in $\FFF$. Therefore there exists $R\in\FFF$ that lies between $P$ and $Q$ and is distinct from these two sets. Then $d(P,R)+d(R,Q)=d(P,Q)$ and both distances $d(P,R)$ and $d(R,Q)$ are less than $n$. By the induction hypothesis, there is a sequence $(R_i)\in\FFF$ such that
$$
P=R_0,\;R=R_j,\;Q=R_n\quad\text{~for some $0<j<n$},
$$
satisfying conditions~(\ref{wg-family}) for $0\leq i<j$ and $j\leq i<n$. It follows that $\FFF$ is a wg-family of sets.
\end{proof}

\vtl
We conclude this section with two examples of infinite partial cubes (more examples are found in~\cite{sO06}).

\begin{example} \label{example Zee}
{\rm Let $\ZZZ$ be the graph on the set $\Zee$ of integers with edges defined by pairs of consecutive integers. 
This graph is a partial cube since its vertex set is isometric to the wg-family of intervals $\{(-\infty,m): m\in\Zee\}$ in $\Zee$.
}
\end{example}

\begin{example} \label{integer lattice example}
{\rm Let us consider $\Zee^n$ as a metric space with respect to the $\ell_1$-metric. The graph $\ZZZ^n$ has $\Zee^n$ as the vertex set; two vertices in $\ZZZ^n$ are connected if they are on the unit distance from each other. We will show in Section~\ref{S:subcubes Cartesian} (Corollary~\ref{product of partial cubes}) that $\ZZZ^n$ is a partial cube.
}
\end{example}

\section{Characterizations} \label{S:characterizations}

Only connected graphs are considered in this section.

\begin{definition}
{\rm Let $G=(V,E)$ be a graph and $d$ be its distance function. For any two adjacent vertices $a,b\in V$ let $W_{ab}$ be the set of vertices that are closer to $a$ than to $b$:
$$
W_{ab}=\{w\in V: d(w,a)<d(w,b)\}.
$$
Following \cite{dE05}, we call the sets $W_{ab}$ and induced subgraphs $\langle W_{ab}\rangle$ {\sl semicubes} of the graph $G$. The semicubes $W_{ab}$ and $W_{ba}$ are called {\sl opposite semicubes}.
}
\end{definition}

\begin{remark}
{\rm The subscript $ab$ in $W_{ab}$ stands for an ordered pair of vertices, not for an edge of $G$. In his original paper~\cite{dD73}, Djokovi\'{c} uses notation $G(a,b)$ (cf.~\cite{mD97}). We use the notation from~\cite{wI00}.
}
\end{remark}

Clearly, two opposite semicubes are disjoint. They can be used to characterize bipartite graphs as follows.

\begin{theorem} \label{bipartite <=> W=W}
A graph $G=(V,E)$ is bipartite if and only if the semicubes $W_{ab}$ and $W_{ba}$ form a partition of $V$ for any edge $ab\in E$.
\end{theorem}

\begin{proof}
Let us recall that a connected graph $G$ is bipartite if and only if for every vertex $x$ there is no edge $ab$ with $d(x,a)=d(x,b)$ (see, for instance,~\cite{aA98}). For any edge $ab\in E$ and vertex $x\in V$ we clearly have
$$
d(x,a)=d(x,b)\quad\Leftrightarrow\quad x\notin W_{ab}\cup W_{ba}.
$$
The result follows.
\end{proof}

The following lemma is instrumental and will be used frequently in the rest of the paper.

\begin{lemma} \label{d(x,b)=d(x,a)+1}
Let $G=(V,E)$ be a graph and $w\in W_{ab}$ for some edge $ab\in E$. Then
$$
d(w,b)=d(w,a)+1.
$$
Accordingly,
$$
W_{ab}=\{w\in V:d(w,b)=d(w,a)+1\}.
$$
\end{lemma}

\begin{proof}
By the triangle inequality, we have
$$
d(w,a)<d(w,b)\leq d(w,a)+d(a,b)=d(w,a)+1.
$$
The result follows, since $d$ takes values in $\Na$.
\end{proof}

There are two binary relations on the set of edges of a graph that play a central role in characterizing partial cubes.

\begin{definition} \label{def theta}
{\rm Let $G=(V,E)$ be a graph and $e=xy$ and $f=uv$ be two edges of $G$.
\roster
	\item[(i)](Djokovi\'{c}~\cite{dD73}) The relation $\hh$ on $E$ is defined by
$$
e\,\hh f\;\eq\;\text{$f$ joins a vertex in $W_{xy}$ with a vertex in $W_{yx}$.}
$$
The notation can be chosen such that $u\in W_{xy}$ and $v\in W_{yx}$.
	\item[(ii)] (Winkler~\cite{pW84}) The relation $\Theta$ on $E$ is defined by
$$
e\,\Theta f\quad\eq\quad d(x,u)+d(y,v)\not=d(x,v)+d(y,u).
$$
\endroster
}
\end{definition}

It is clear that both relations $\hh$ and $\Theta$ are reflexive and $\Theta$ is symmetric.

\begin{lemma}
The relation $\hh$ is a symmetric relation on $E$.
\end{lemma}

\begin{proof}
Suppose that $xy\,\hh\,uv$ with $u\in W_{xy}$ and $v\in W_{yx}$. By Lemma~\ref{d(x,b)=d(x,a)+1} and the triangle inequality, we have
\begin{align*}
d(u,x)&=d(u,y)-1\leq d(u,v)+d(v,y)-1=d(v,y)= \\
&=d(v,x)-1\leq d(v,u)+d(u,x)-1=d(u,x).
\end{align*}
Hence, $d(u,x)=d(v,x)-1$ and $d(v,y)=d(u,y)-1$. Therefore, $x\in W_{uv}$ and $y\in W_{vu}$. It follows that $uv\,\hh\,xy$.
\end{proof}

\begin{lemma} \label{theta in Theta}
$\hh\SB\Theta$.
\end{lemma}

\begin{proof}
Suppose that $xy\,\hh\,uv$ with $u\in W_{xy},\;v\in W_{yx}$. By Lemma~\ref{d(x,b)=d(x,a)+1},
$$
d(x,u)+d(y,v)=d(x,v)-1+d(y,u)-1\not=d(x,v)+d(y,u).
$$
Hence, $xy\,\Theta\,uv$.
\end{proof}

\begin{example}
{\rm It is easy to verify that $\hh$ is the identity relation on the set of edges of the cycle $C_3$. On the other hand, any two edges of $C_3$ stand in the relation $\Theta$. Thus, $\hh\not=\Theta$ in this case.
}
\end{example}

Bipartite graphs can be characterized in terms of relations $\hh$ and $\Theta$ as follows.

\begin{theorem} \label{bipartite <=> theta=Theta}
A graph $G=(V,E)$ is bipartite if and only if $\hh=\Theta$.
\end{theorem}

\begin{proof}
(Necessity.) Suppose that $G$ is a bipartite graph, two edges $xy$ and $uv$ stand in the relation $\Theta$, that is,
$$
d(x,u)+d(y,v)\not=d(x,v)+d(y,u),
$$
and that edges $xy$ and $uv$ do not stand in the relation $\hh$. By Theorem~\ref{bipartite <=> W=W}, we may assume that $u,v\in W_{xy}$. By Lemma~\ref{d(x,b)=d(x,a)+1}, we have
$$
d(x,u)+d(y,v)=d(y,u)-1+d(x,v)+1=d(x,v)+d(y,u),
$$
a contradiction. It follows that $\Theta\SB\hh$. By Lemma~\ref{theta in Theta}, $\hh=\Theta$.

\vtl
(Sufficiency.) Suppose that $G$ is not bipartite. By Theorem~\ref{bipartite <=> W=W}, there is an edge $xy$ such that $W_{xy}\cup W_{yx}$ is a proper subset of $V$. Since $G$ is connected, there is an edge $uv$ with $u\notin W_{xy}\cup W_{yx}$ and $v\in W_{xy}\cup W_{yx}$. Clearly, $uv$ does not stand in the relation $\hh$ to $xy$. On the other hand,
$$
d(x,u)+d(y,v)\not=d(x,v)+d(y,u),
$$
since $u\notin W_{xy}\cup W_{yx}$ and $v\in W_{xy}\cup W_{yx}$. Thus, $xy\,\Theta\,uv$, a contradiction, since we assumed that $\hh=\Theta$.
\end{proof}

\vtl
By Theorem~\ref{bipartite <=> theta=Theta}, the relations $\hh$ and $\Theta$ coincide on bipartite graphs. For this reason we use the relation $\hh$ in the rest of the paper.

\vtl
\begin{lemma} \label{convexity => partitions}
Let $G=(V,E)$ be a bipartite graph such that all its semicubes are convex sets. Then two edges $xy$ and $uv$ stand in the relation $\hh$ if and only if the corresponding pairs of mutually opposite semicubes form equal partitions of $V$:
$$
xy\,\hh\,uv\quad\eq\quad\{W_{xy},W_{yx}\}=\{W_{uv},W_{vu}\}.
$$
\end{lemma}

\begin{proof}
(Necessity) We assume that the notation is chosen such that $u\in W_{xy}$ and $v\in W_{yx}$. Let $z\in W_{xy}\cap W_{vu}$. By Lemma~\ref{d(x,b)=d(x,a)+1}, $d(z,u)=d(z,v)+d(v,u)$. Since $z,u\in W_{xy}$ and $W_{xy}$ is convex, we have $v\in W_{xy}$, a contradiction to the assumption that $v\in W_{yx}$. Thus $W_{xy}\cap W_{vu}=\es$. Since two opposite semicubes in a bipartite graph form a partition of $V$, we have $W_{uv}=W_{xy}$ and $W_{vu}=W_{yx}$.

A similar argument shows that $W_{uv}=W_{yx}$ and $W_{vu}=W_{xy}$, if $u\in W_{yx}$ and $v\in W_{xy}$.

\vtl
(Sufficiency.) Follows from the definition of the relation $\hh$.
\end{proof}

\vtl
We need another general property of the relation $\hh$ (cf.~Lemma~2.2 in~\cite{wI00}).

\begin{lemma} \label{edges in path}
Let $P$ be a shortest path in a graph $G$. Then no two distinct edges of $P$ stand in the relation $\hh$.
\end{lemma}

\begin{proof}
Let $i<j$ and $x_i x_{i+1}$ and $x_j x_{j+1}$ be two edges in a shortest path $P$ from $x_0$ to $x_n$. Then 
$$
d(x_i,x_j)<d(x_i,x_{j+1})\quad\text{and}\quad d(x_{i+1},x_j)<d(x_{i+1},x_{j+1}),
$$
so $x_i,x_{i+1}\in W_{x_jx_{j+1}}$. It follows that edges $x_i x_{i+1}$ and $x_j x_{j+1}$ do not stand in the relation $\hh$.
\end{proof}

The converse statement is true for bipartite graphs (we omit the proof); a counterexample is the cycle $C_5$ which is not bipartite.

\begin{lemma} \label{convexity <=> transitivity}
Let $G=(V,E)$ be a bipartite graph. The following statements are equivalent
\roster
	\item[{\rm(i)}] All semicubes of $G$ are convex.
	\item[{\rm(ii)}] The relation $\hh$ is an equivalence relation on $E$.
\endroster
\end{lemma}

\begin{proof}
(i) $\Rightarrow$ (ii). Follows from Lemma~\ref{convexity => partitions}.

\vtl
(ii) $\Rightarrow$ (i). Suppose that $\hh$ is transitive and there is a nonconvex semicube $W_{ab}$. Then there are two vertices $u,v\in W_{ab}$ and a shortest path $P$ from $u$ to $v$ that intersects $W_{ba}$. This path contains two distinct edges $e$ and $f$ joining vertices of semicubes $W_{ab}$ and $W_{ba}$. The edges $e$ and $f$ stand in the relation $\hh$ to the edge $ab$. By transitivity of $\hh$, we have $e\,\hh f$. This contradicts the result of Lemma~\ref{edges in path}. Thus all semicubes of $G$ are convex.
\end{proof}

\vtl
We now establish some basic properties of partial cubes.

\vtl
\begin{theorem} \label{partial cube => properties}
Let $G=(V,E)$ be a partial cube. Then
\roster
	\item[{\rm(i)}] $G$ is a bipartite graph.
	\item[{\rm(ii)}] Each pair of opposite semicubes form a partition of $V$.
	\item[{\rm(iii)}] All semicubes are convex subsets of $V$.
	\item[{\rm(iv)}] $\hh$ is an equivalence relation on $E$.
\endroster
\end{theorem}

\begin{proof}
We may assume that $G$ is an isometric subgraph of some hypercube $\HHH(X)$, that is, $G=(\FFF,E_\FFF)$ for a wg-family $\FFF$ of finite subsets of $X$.
\vtl
(i) It suffices to note that if two sets in $\HHH(X)$ are connected by an edge then they have different parity. Thus, $\HHH(X)$ is a bipartite graph and so is $G$.
\vtl
(ii) Follows from (i) and Theorem~\ref{bipartite <=> W=W}.
\vtl
(iii) Let $W_{AB}$ be a semicube of $G$. By Lemma~\ref{d(x,b)=d(x,a)+1} and Theorem~\ref{lattice=metric betweenness}, we have
$$
W_{AB}=\{S\in\FFF:S\cap B\SB A\SB S\cup B\}.
$$
Let $Q,R\in W_{AB}$ and $P$ be a vertex of $G$ such that
$$
d(Q,P)+d(P,R)=d(Q,R).
$$
By Theorem~\ref{lattice=metric betweenness},
$$
Q\cap R\SB P\SB Q\cup R.
$$
Since $Q,R\in W_{AB}$, we have
$$
Q\cap B\SB A\SB Q\cup B\quad\text{and}\quad R\cap B\SB A\SB R\cup B,
$$
which implies
$$
P\cap B\SB(Q\cup R)\cap B\SB A\SB(Q\cap R)\cup B\SB S\cup B.
$$
Hence, $P\in W_{AB}$, and the result follows.

\vtl
(iv) Follows from (iii) and Lemma~\ref{convexity <=> transitivity}.
\end{proof}

\begin{remark} \label{half-space}
{\rm Since semicubes of a partial cube $G=(V,E)$ are convex subsets of the metric space $V$, they are {\sl half-spaces} in $V$~\cite{mV93}. This terminology is used in~\cite{vC88,vC94}.
}
\end{remark}

The following theorem presents four characterizations of partial cubes. The first two are due to Djokovi\'{c}~\cite{dD73} and Winkler~\cite{pW84} (cf. Theorem~2.10 in~\cite{wI00}).

\begin{theorem} \label{DWT}
Let $G=(V,E)$ be a connected graph. The following statements are equivalent:
\roster
	\item[{\rm(i)}] $G$ is a partial cube.
	\item[{\rm(ii)}] $G$ is bipartite and all semicubes of $G$ are convex.
	\item[{\rm(iii)}] $G$ is bipartite and $\hh$ is an equivalence relation.
	\item[{\rm(iv)}] $G$ is bipartite and, for all $xy,uv\in E$, 
\begeq \label{W=W}
xy\,\hh\,uv\quad\imp\quad\{W_{xy},W_{yx}\}=\{W_{uv},W_{vu}\}.
\edeq
	\item[{\rm(v)}] $G$ is bipartite and, for any pair of adjacent vertices of $G$, there is a unique pair of opposite semicubes separating these two vertices.
\endroster
\end{theorem}

\begin{proof}
By Lemma~\ref{convexity <=> transitivity}, the statements (ii) and (iii) are equivalent and, by Theorem~\ref{partial cube => properties}, (i) implies both (ii) and (iii).
\vtl

(iii) $\Rightarrow$ (i). By Theorem~\ref{bipartite <=> W=W}, each pair $\{W_{ab},W_{ba}\}$ of opposite semicubes of $G$ form a partition of $V$. We orient these partitions by calling, in an arbitrary way, one of the two opposite semicubes in each partition a {\sl positive semicube}. Let us assign to each $x\in V$ the set $W^+(x)$ of all positive semicubes containing $x$. In the next paragraph we prove that the family $\FFF=\{W^+(x)\}_{x\in V}$ is well graded and that the assignment $x\mapsto W^+(x)$ is an isometry between $V$ and $\FFF$.

\vtl
Let $x$ and $y$ be two distinct vertices of $G$. We say that a positive semicube $W_{ab}$ {\sl separates} $x$ and $y$ if either $x\in W_{ab},\;y\in W_{ba}$ or $x\in W_{ba},\;y\in W_{ab}$. It is clear that $W_{ab}$ separates $x$ and $Y$ if and only if $W_{ab}\in W^+(x)\Delta W^+(y)$. Let $P$ be a shortest path $x_0=x,x_1,\ldots,x_n=y$ from $x$ to $y$. By Lemma~\ref{edges in path}, no two distinct edges of $P$ stand in the relation $\hh$. By Lemma~\ref{convexity => partitions}, distinct edges of $P$ define distinct positive semicubes; clearly, these semicubes separate $x$ and $y$. Let $W_{ab}$ be a positive semicube separating $x$ and $y$, and, say, $x\in W_{ab}$ and $y\in W_{ba}$. There is an edge $f\in P$ that joins vertices in $W_{ab}$ and $W_{ba}$. Hence, $f$ stands in the relation $\hh$ to $ab$ and, by Lemma~\ref{convexity => partitions}, $W_{ab}$ is defined by $f$. It follows that any semicube in $W^+(x)\Delta W^+(y)$ is defined by a unique edge in $P$ and any edge in $P$ defines a semicube in $W^+(x)\Delta W^+(y)$. Therefore, $d(W^+(x),W^+(y))=d(x,y)$, that is $x\mapsto W^+(x)$ is an isometry. Clearly, $\FFF$ is a wg-family of sets.

By Theorem~\ref{wg=finite wg}, the family $\FFF$ is isometric to a wg-family of finite sets. Hence, $G$ is a partial cube.

\vtl
(iv) $\Rightarrow $ (ii). Suppose that there exist an edge $ab$ such that semicube $W_{ba}$ is not convex. Let $p$ and $q$ be two vertices in $W_{ba}$ such that there is a shortest path $P$ from $p$ to $q$ that intersects $W_{ab}$. There are two distinct edges $xy$ and $uv$ in $P$ such that $x,u\in W_{ab}$ and $y,v\in W_{ba}$. Since $ab\,\hh\,xy$ and $ab\,\hh\,uv$, we have, by~(\ref{W=W}),
$$
W_{ab}=W_{xy}=W_{uv}.
$$
Hence, $u\in W_{xy}$ and $v\in W_{yx}$. By Lemma~\ref{d(x,b)=d(x,a)+1},
$$
d(x,u)=d(x,v)-1=1+d(v,y)-1=d(v,y),
$$
a contradiction, since $P$ is a shortest path from $p$ to $q$.

\vtl
(ii) $\Rightarrow $ (iv). Follows from Lemma~\ref{convexity => partitions}.

\vtl
It is clear that (iv) and (v) are equivalent.
\end{proof}

\section{Fundamental sets in partial cubes} \label{S:fundamental sets}

Semicubes played an important role in the previous section. In this section we introduce three more classes of useful subsets of graphs. We also establish one more characterization of partial cubes.

\vtl
Let $G=(V,E)$ be a connected graph. For a given edge $e=ab\in E$, we define the following sets (cf.~\cite{wI00,hM80}):
\begin{align*}
F_{ab} &= \{f\in E : e\,\hh f\} = \{uv\in E: u\in W_{ab},v\in W_{ba}\}, \\
U_{ab} &= \{w\in W_{ab}:\text{$w$ is adjacent to a vertex in $W_{ba}$}\}, \\
U_{ba} &= \{w\in W_{ba}:\text{$w$ is adjacent to a vertex in $W_{ab}$}\}.
\end{align*}
The five sets are schematically shown in Figure~\ref{fundamental sets}.

{\begin{figure}[h!]
\vspace{0.3in}
\centerline{\includegraphics{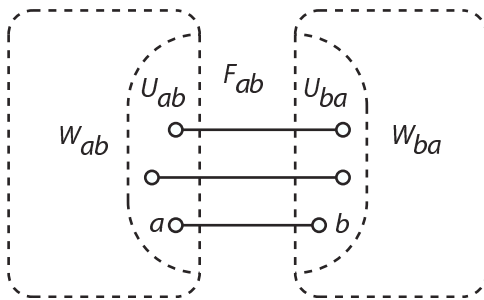}}
\vspace{0.3in}
\caption {Fundamental sets in a partial cube.} \label{fundamental sets} 
\end{figure}
}

\begin{remark}
{\rm In the case of a partial cube $G=(V,E)$, the semicubes $W_{ab}$ and $W_{ba}$ are complementary half-spaces in the metric space $V$ (cf. Remark~\ref{half-space}). Then the set $F_{ab}$ can be regarded as a `hyperplane' separating these half-spaces (see~\cite{sO06} where this analogy is formalized in the context of hyperplane arrangements).
}
\end{remark}

The following theorem generalizes the result obtained in~\cite{hM80} for median graphs (see also~\cite{wI00}).

\begin{theorem} \label{theorem matching=isomorphism}
Let $ab$ be an edge of a connected bipartite graph $G$. If the semicubes $W_{ab}$ and $W_{ba}$ are convex, then the set $F_{ab}$ is a matching and induces an isomorphism between the graphs $\langle U_{ab}\rangle$ and $\langle U_{ba}\rangle$.
\end{theorem}

\begin{proof}
Suppose that $F_{ab}$ is not a matching. Then there are distinct edges $xu$ and $xv$ with, say, $x\in U_{ab}$ and $u,v\in U_{ba}$. By the triangle inequality, $d(u,v)\leq 2$. Since $G$ does not have triangles, $d(u,v)\not=1$. Hence, $d(u,v)=2$, which implies that $x$ lies between $u$ and $v$. This contradicts convexity of $W_{ba}$, since $x\in W_{ab}$. Therefore $F_{ab}$ is a matching.

To show that $F_{ab}$ induces an isomorphism, let $xy,uv\in F_{ab}$ and $xu\in E$, where $x,u\in U_{ab}$ and $y,v\in U_{ba}$. Since $G$ does not have odd cycles, $d(v,y)\not=2$. By the triangle inequality,
$$
d(v,y)\leq d(v,u)+d(u,x)+d(x,y)=3.
$$
Since $W_{ba}$ is convex, $d(v,y)\not=3$. Thus $d(v,y)=1$, that is, $vy$ is an edge. The result follows by symmetry.
\end{proof}

By Theorem~\ref{DWT}(ii), we have the following corollary.

\begin{corollary} \label{matching=isomorphism}
Let $G=(V,E)$ be a partial cube. For any edge $ab$ the set $F_{ab}$ is a matching and induces an isomorphism between induced graphs $\langle U_{ab}\rangle$ and $\langle U_{ba}\rangle$.
\end{corollary}

{\begin{figure}[h!]
\centerline{\includegraphics{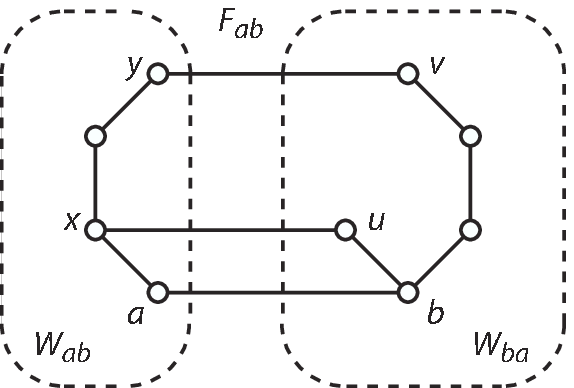}}
\caption {Graph $G$.} \label{example3} 
\end{figure}
}

\begin{example}
{\rm Let $G$ be the graph depicted in Figure~\ref{example3}. The set $$F_{ab}=\{ab,xu,yv\}$$ is a matching and defines an isomorphism between the graphs induced by subsets $U_{ab}=\{a,x,y\}$ and $U_{ba}=\{b,u,v\}$. The set $W_{ba}$ is not convex, so $G$ is not a partial cube. Thus the converse of Corollary~\ref{matching=isomorphism} does not hold.
}
\end{example}

We now establish another characterization of partial cubes that utilizes a geometric property of families $F_{ab}$.

\begin{theorem} \label{F-characterization}
For a connected graph $G$ the following statements are equivalent:
\roster
	\item[{\rm(i)}] $G$ is a partial cube.
	\item[{\rm (ii)}] $G$ is bipartite and
\begeq \label{rectangle}
d(x,u)=d(y,v)\quad\text{and}\quad d(x,v)=d(y,u),
\edeq
for any $ab\in E$ and $xy,uv\in F_{ab}$.
\endroster
\end{theorem}

\begin{proof}
(i)$\Rightarrow$(ii). We may assume that $x,u\in W_{ab}$ and $y,v\in W_{ba}$. Since $\hh$ is an equivalence relation, we have $xy\,\hh\,uv\,\hh ab$. By Lemma~\ref{convexity => partitions}, $W_{uv}=W_{xy}=W_{ab}$. By Lemma~\ref{d(x,b)=d(x,a)+1},
$$
d(x,u)=d(x,v)-1=d(v,y)+1-1=d(y,v).
$$
We also have
$$
d(x,v)=d(y,v)+1=d(y,u),
$$
by the same lemma.

\vtl
(ii)$\Rightarrow$(i). Suppose that $G$ is not a partial cube. Then, by Theorem~\ref{DWT}, there exist an edge $ab$ such that, say, semicube $W_{ba}$ is not convex. Let $p$ and $q$ be two vertices in $W_{ba}$ such that there is a shortest path $P$ from $p$ to $q$ that intersects $W_{ab}$. Let $uv$ be the first edge in $P$ which belongs to $F_{ab}$ and $xy$ be the last edge in $P$ with the same property (see Figure~\ref{example5}).

{\begin{figure}[h!]
\centerline{\includegraphics{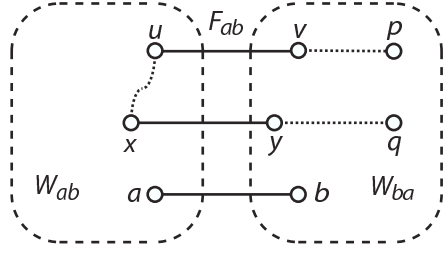}}
\vspace{0.1in}
\caption {An illustration to the proof of theorem~\ref{F-characterization}.} 
\label{example5} 
\end{figure}
}
\noindent
Since $P$ is a shortest path, we have
$$
d(v,y)=d(v,u)+d(u,x)+d(x,y)\not=d(x,u),
$$
which contradicts condition~(\ref{rectangle}). Thus all semicubes of $G$ are convex. By Theorem~\ref{DWT}, $G$ is a partial cube.
\end{proof}

\vtl
\begin{remark}
{\rm One can say that four vertices satisfying conditions~(\ref{rectangle}) define a {\sl rectangle} in $G$. Then Theorem~\ref{F-characterization} states that a connected graph is a partial cube if and only if it is bipartite and for any edge $ab$ pairs of edges in $F_{ab}$ define rectangles in $G$. 
}
\end{remark}

\section{Dimensions of partial cubes} \label{S:dimensions}

There are many different ways in which a given partial cube can be isometrically embedded into a hypercube. For instance, the graph $K_2$ can be isometrically embedded in different ways into any hypercube $\HHH(X)$ with $|X|>2$.

\vtl
Following Djokovi\'{c}~\cite{dD73} (see also~\cite{mD97}), we define the {\sl isometric dimension}, $\dim_I(G)$, of a partial cube $G$ as the minimum possible dimension of a hypercube $\HHH(X)$ in which $G$ is isometrically embeddable. Recall (see Section~\ref{S:hypercubes}) that the dimension of $\HHH(X)$ is the cardinality of the set $X$.

\begin{theorem} \label{Djokovic dimension}
{\rm(Theorem~2 in~\cite{dD73}.)} Let $G=(V,E)$ be a partial cube. Then
\begeq \label{isometric dimension}
\dim_I(G)=|E\slash\hh|,
\edeq
where $\hh$ is Djokovi\'{c}'s equivalence relation on $E$ and $E\slash\hh$ is the set of its equivalence classes {\rm(}the quotient-set{\rm)}.
\end{theorem}

The quotient-set $E\slash\hh$ can be identified with the family of all distinct sets $F_{ab}$ (see Section~\ref{S:fundamental sets}). If $G$ is a finite partial cube, we may consider it as an isometric subgraph of some hypercube $Q_n$. Then the edges in each family $F_{ab}$ are parallel edges in $Q_n$ (cf. Theorem~\ref{F-characterization}). This observation essentially proves~(\ref{isometric dimension}) in the finite case.

\vtl
Let $G$ be a partial cube on a set $X$. The vertex set of $G$ is a wg-family $\FFF$ of finite subsets of $X$ (see Section~\ref{S:hypercubes}). We define the {\sl retraction} of $\FFF$ as a family $\FFF'$ of subsets of $X'=\cup\,\FFF\setminus\cap\,\FFF$ consisting of the intersections of sets in $\FFF$ with $X'$. It is clear that $\FFF'$ satisfies conditions
\begeq \label{cap-cup}
\cap\,\FFF'=\es\quad\text{and}\quad\cup\FFF'=X'.
\edeq

\begin{proposition} \label{retraction}
The partial cubes induced by a wg-family $\FFF$ and its retraction $\FFF'$ are isomorphic.
\end{proposition}

\begin{proof}
It suffices to prove that metric spaces $\FFF$ and $\FFF'$ are isometric. Clearly, $\aa:P\mapsto P\cap X'$ is a mapping from $\FFF$ onto $\FFF'$. For $P,Q\in\FFF$, we have
$$
(P\cap X')\Delta(Q\cap X')=(P\Delta Q)\cap X'=(P\Delta Q)\cap(\cup\FFF\setminus\cap\FFF)=P\Delta Q.
$$
Thus, $d(\aa(P),\aa(Q))=d(P,Q)$. Consequently, $\aa$ is an isometry.
\end{proof}

\vtl
Let $G$ be a partial cube on some set $X$ induced by a wg-family $\FFF$ satisfying conditions~(\ref{cap-cup}), and let $PQ$ be an edge of $G$. By definition, there is $x\in X$ such that $P\Delta Q=\{x\}$. The following two lemmas are instrumental.

\begin{lemma} \label{points=semicubes}
Let $PQ$ be an edge of a partial cube $G$ on $X$ and let $P\Delta Q=\{x\}$. The two sets
$$
\{R\in\FFF: x\in R\}\quad\text{and}\quad\{R\in\FFF: x\notin R\}
$$
form the same bipartition of the family $\FFF$ as semicubes $W_{PQ}$ and $W_{QP}$.
\end{lemma}

\begin{proof}
We may assume that $Q=P+\{x\}$. Then, for any $R\in\FFF$,
$$
R\Delta Q = R\Delta(P+\{x\}) = \begin{cases}
	(R\Delta P)+\{x\}, &\text{if $x\in R$,}\\
	R\Delta P, &\text{if $x\notin R$.}
\end{cases}
$$
Hence, $|R\Delta P|<|R\Delta Q|$ if and only if $x\in R$. It follows that
$$
W_{PQ}=\{R\in\FFF: x\in R\}.
$$
A similar argument shows that $W_{QP}=\{R\in\FFF: x\notin R\}$.
\end{proof}

\vtl
\begin{lemma} \label{semicubes=points}
If $\FFF$ is a wg-family of sets satisfying conditions~{\rm(\ref{cap-cup})}, then for any $x\in X$ there are sets $P,Q\in\FFF$ such that $P\Delta Q=\{x\}$.
\end{lemma}

\begin{proof}
By conditions~\ref{cap-cup}, for a given $x\in X$ there are sets $S$ and $T$ in $\FFF$ such that $x\in S$ and $x\notin T$. Let $R_0=S,R_1,\ldots,R_n=T$ be a sequence of sets in $\FFF$ satisfying conditions~(\ref{wg-family}). It is clear that there is $i$ such that $x\in R_i$ and $x\notin R_{i+1}$. Hence, $R_i\Delta R_{i+1}=\{x\}$, so we can choose $P=R_i$ and $Q=R_{i+1}$.
\end{proof}

\vtl
By Lemmas~\ref{points=semicubes} and~\ref{semicubes=points}, there is one-to-one correspondence between the set $X$ and the quotient-set $E\slash\hh$. From Theorem~\ref{Djokovic dimension} we obtain the following result.

\vtl
\begin{theorem} \label{X-dimension}
Let $\FFF$ be a wg-family of finite subsets of a set $X$ such that $\cap\,\FFF=\es$ and $\cup\,\FFF=X$, and let $G$ be a partial cube on $X$ induced by $\FFF$. Then
$$
\dim_I(G)=|X|.
$$
\end{theorem}

\vtl
Clearly, a graph which is isometrically embeddable into a partial cube is a partial cube itself. We will show in Section~\ref{S:subcubes Cartesian} (Corollary~\ref{product of partial cubes}) that the integer lattice $\ZZZ^n$ is a partial cube. Thus a graph which is isometrically embeddable into an integer lattice is a partial cube. It follows that a finite graph is a partial cube if and only if it is embeddable in some integer lattice. Examples of infinite partial cubes isometrically embeddable into a finite dimensional integer lattice are found in~\cite{sO06}. 

\vtl
We call the minimum possible dimension $n$ of an integer lattice $\ZZZ^n$, in which a given graph $G$ is isometrically embeddable, its {\sl lattice dimension} and denote it $\dim_Z(G)$. The lattice dimension of a partial cube can be expressed in terms of maximum matchings in so-called semicube graphs~\cite{dE05}.

\begin{definition}
{\rm The {\sl semicube graph} $\text{Sc}(G)$ has all semicubes in $G$ as the set of its vertices. Two vertices $W_{ab}$ and $W_{cd}$ are connected in $\text{Sc}(G)$ if
\begin{equation} \label{compatible1}
W_{ab}\cup W_{cd}=V\quad\text{and}\quad W_{ab}\cap W_{cd}\not=\emptyset.
\end{equation}
}
\end{definition}

If $G$ is a partial cube, then condition~(\ref{compatible1}) is equivalent to each of the two equivalent conditions:
\begin{equation} \label{compatible2}
W_{ba}\subset W_{cd}\quad\Leftrightarrow\quad W_{dc}\subset W_{ab},
\end{equation}
where $\subset$ stands for the proper inclusion.

\vtl
\begin{theorem} \label{lattice=isometric-|M|}
{\rm (Theorem~1 in~\cite{dE05}.)} Let $G$ be a finite partial cube. Then
$$
\dim_Z(G) = \dim_I(G)-|M|,
$$
where $M$ is a maximum matching in the semicube graph $\text{{\rm Sc}}(G)$.
\end{theorem}

\begin{example}
{\rm Let $G$ be the graph shown in Figure~\ref{example0}. It is easy to see that
$$
\dim_I(G)=3\quad\text{and}\quad\dim_Z(G)=2.
$$
}
\end{example}

\begin{example} \label{tree dimensions}
{\rm Let $T$ be a tree with $n$ edges and $m$ leaves. Then
$$
\dim_I(T)=n\quad\text{and}\quad\dim_Z(T)=\lceil m\slash 2\rceil
$$
(cf.~\cite{mD97} and~\cite{fH78}, respectively).
}
\end{example}

\begin{example}
{\rm For the cycle $C_6$ we have (see Figure~\ref{example15})
$$
\dim_I(C_6)=\dim_Z(C_6)=3.
$$
}
\end{example}

\section{Subcubes and Cartesian products} \label{S:subcubes Cartesian}

Let $G$ be a partial cube. We say that $G'$ is a {\sl subcube} of $G$ if it is an isometric subgraph of $G$. 

\vtl
Clearly, a subcube is itself a partial cube. The converse does not hold; a subgraph of a graph $G$ can be a partial cube but not an isometric subgraph of $G$ (cf. Example~\ref{nonisometric path}).

\vtl
If $G'$ is a subcube of a partial cube $G$, then $\dim_I(G')\leq\dim_I(G)$ and $\dim_Z(G')\leq\dim_Z(G)$. In general, the two inequalities are not strict. For instance, the cycle $C_6$ is an isometric subgraph of the cube $Q_3$ (see Figure~\ref{example15}) and
$$
\dim_I(C_6)=\dim_Z(C_6)=\dim_I(Q_3)=\dim_Z(Q_3)=3.
$$

Semicubes of a partial cube are examples of subcubes. Indeed, by Theorem~\ref{DWT}, semicubes are convex subgraphs and therefore isometric. In general, the converse is not true; a path connecting two opposite vertices in $C_6$ is an isometric subgraph but not a convex one.

\vtl
Another common way of constructing new partial cubes from old ones is by forming their Cartesian products (see~\cite{wI00} for details and proofs).

\begin{definition} \label{product}
{\rm Given two graphs $G_1=(V_1,E_1)$ and $G_2=(V_2,E_2)$, their {\sl Cartesian product} 
$$
G=G_1\Box\,G_2
$$
has vertex set $V=V_1\times V_2$; a vertex $u=(u_1,u_2)$ is adjacent to a vertex $v=(v_1,v_2)$ if and only if $u_1v_1\in E_1$ and $u_2=v_2$, or $u_1=v_1$ and $u_2v_2\in E_2$.
}
\end{definition}

The operation $\Box$ is associative, so we can write
$$
G=G_1\Box\cdots\Box\,G_n=\prod_{i=1}^n G_i
$$
for the Cartesian product of graphs $G_1,\ldots,G_n$. A Cartesian product $\prod_{i=1}^n G_i$ is connected if and only if the factors are connected. Then we have

\begeq \label{sum of distances}
d_G(u,v)=\sum_{i=1}^n d_{G_i}(u_i,v_i).
\edeq

\begin{example} \label{hypercube product}
{\rm Let $\{X_i\}_{i=1}^n$ be a family of sets and $Y=\sum_{i=1}^n$ be their sum. Then the Cartesian product of the hypercubes $\HHH(X_i)$ is isomorphic to the hypercube $\HHH(Y)$. The isomorphism is established by the mapping 
$$
f:(P_1,\ldots,P_n)\mapsto \sum_{i=1}^n P_i.
$$
}
\end{example}

Formula~(\ref{sum of distances}) yields immediately the following results.

\begin{proposition}
Let $H_i$ be isometric subgraphs of graphs $G_i$ for all $1\leq i\leq n$. Then the Cartesian product $\prod_{i=1}^n H_i$ is an isometric subgraph of the Cartesian product $\prod_{i=1}^n G_i$.
\end{proposition}

\begin{corollary} \label{product of partial cubes}
The Cartesian product of a finite family of partial cubes is a partial cube. In particular, the integer lattice $\ZZZ^n$ (cf. Examples~\ref{example Zee} and~\ref{integer lattice example}) is a partial cube.
\end{corollary}

The results of the next two theorems can be easily extended to arbitrary finite products of finite partial cubes.

\begin{theorem} \label{Cartesian dimension}
Let $G=G_1\Box\,G_2$ be the Cartesian product of two finite partial cubes. Then
$$
\dim_I(G)=\dim_I(G_1)+\dim_I(G_2).
$$
\end{theorem}

\begin{proof}
We may assume that $G_1$ (resp. $G_2$) is induced by a wg-family $\FFF_1$ (resp. $\FFF_2$) of subsets of a finite set $X_1$ (resp. $X_2$) such that $\cap\,\FFF_1=\es$ and $\cup\,\FFF_1= X_1$ (resp. $\cap\,\FFF_2=\es$ and $\cup\,\FFF_2= X_1$) (see Section~\ref{S:dimensions}). By Theorem~\ref{X-dimension},
$$
\dim_I(G_1)=|X_1|\quad\text{and}\quad\dim_I(G_2)=|X_2|.
$$
It is clear that the graph $G$ is induced by the wg-family $\FFF=\FFF_1+\FFF_2$ of subsets of the set $X=X_1+X_2$ (cf. Example~\ref{hypercube product}) with $\cap\,\FFF=\es,\;\cup\,\FFF= X$. By Theorem~\ref{X-dimension}, 
$$
\dim_I(G)=|X|=|X_1|+|X_2|=\dim_I(G_1)+\dim_I(G_2).
$$
\end{proof}

\begin{theorem}
Let $G=(V,E)$ be the Cartesian product of two finite partial cubes $G_1=(V_1,E_1)$ and $G_2=(V_2,E_2)$. Then
$$
\dim_Z(G)=\dim_Z(G_1)+\dim_Z(G_2).
$$
\end{theorem}

\begin{proof}
Let $W_{(a,b)(c,d)}$ be a semicube of the graph $G$. There are two possible cases:

\vtl
(i) $c=a,\;bd\in E_2$. Let $(x,y)$ be a vertex of G. Then, by~(\ref{sum of distances}),
$$
d_G((x,y),(a,b))=d_{G_1}(x,a)+d_{G_2}(y,b)
$$
and
$$
d_G((x,y),(c,d))=d_{G_1}(x,c)+d_{G_2}(y,d).
$$
Hence,
$$
d_G((x,y),(a,b))<d_G((x,y),(c,d))\quad\Leftrightarrow\quad d_{G_2}(y,b)<d_{G_2}(y,d).
$$
It follows that
\begeq \label{product 1}
W_{(a,b)(c,d)}=V_1\times W_{bd}.
\edeq

\vtl
(ii) $d=b,\;ac\in E_1$. Like in (i), we have
\begeq \label{product 2}
W_{(a,b)(c,d)}=W_{ac}\times V_2.
\edeq

\vtl
Clearly, two semicubes given by~(\ref{product 1}) form an edge in the semicube graph $\text{Sc}(G)$ if and only if their second factors form an edge in the semicube graph $\text{Sc}(G_2)$. The same is true for semicubes in the form~(\ref{product 2}) with respect to their first factors. It is also clear that semicubes in the form~(\ref{product 1}) and in the form~(\ref{product 2}) are not connected by an edge in $\text{Sc}(G)$. Therefore the semicube graph $\text{Sc}(G)$ is isomorphic to the disjoint union of semicube graphs $\text{Sc}(G_1)$ and $\text{Sc}(G_2)$. If $M_1$ is a maximum matching in $\text{Sc}(G_1)$ and $M_2$ is a maximum matching in $\text{Sc}(G_2)$, then $M=M_1\cup M_2$ is a maximum matching in $\text{Sc}(G)$. The result follows from theorems~\ref{lattice=isometric-|M|} and~\ref{Cartesian dimension}.
\end{proof}

\begin{remark}
{\rm The result of Corollary~\ref{product of partial cubes} does not hold for infinite Cartesian products of partial cubes, as these products are disconnected. On the other hand, it can be shown that arbitrary {\sl weak} Cartesian products (connected components of Cartesian products~\cite{wI00}) of partial cubes are partial cubes.
}
\end{remark}

\section{Pasting partial cubes} \label{S:pasting}

In this section we use the set pasting technique~\cite[ch.I, \S2.5]{nB66} to build new partial cubes from old ones.

\vtl
Let $G_1=(V_1,E_1)$ and $G_2=(V_2,E_2)$ be two graphs, $H_1=(U_1,F_1)$ and $H_2=(U_2,F_2)$ be two isomorphic subgraphs of $G_1$ and $G_2$, respectively, and $\psi:U_1\rightarrow U_2$ be a bijection defining an isomorphism between $H_1$ and $H_2$. The bijection $\psi$ defines an equivalence relation $R$ on the sum $V_1+V_2$ as follows: any element in $(V_1\setminus U_1)\cup(V_2\setminus U_2)$ is equivalent to itself only and elements $u_1\in U_1$ and $u_2\in U_2$ are equivalent if and only if $u_2=\psi(u_1)$. We say that the quotient set $V=(V_1+V_2)\slash R$ is obtained by {\sl pasting together the sets $V_1$ and $V_2$ along the subsets $U_1$ and $U_2$}. Since the graphs $H_1$ and $H_2$ are isomorphic, the pasting of the sets $V_1$ and $V_2$ can be naturally extended to a pasting of sets of edges $E_1$ and $E_2$ resulting in the set $E$ of edges joining vertices in $V$. We say that the graph $G=(E,V)$ is obtained by {\sl pasting together the graphs $G_1$ and $G_2$ along the isomorphic subgraphs $H_1$ and $H_2$}. The pasting construction allows for identifying in a natural way the graphs $G_1$ and $G_2$ with subgraphs of $G$, and the isomorphic graphs $H_1$ and $H_2$ with a common subgraph $H$ of both graphs $G_1$ and $G_2$. We often follow this convention below.

\begin{remark}
{\rm Note that in the above construction the resulting graph $G$ depends not only on graphs $G_1$ and $G_2$ and their isomorphic subgraphs $H_1$ and $H_2$ but also on the bijection $\psi$ defining an isomorphism from $H_1$ onto $H_2$ (see the drawings in Figures~\ref{example6} and~\ref{example7}). 
}
\end{remark}

{\begin{figure}[h!]
\centerline{\includegraphics{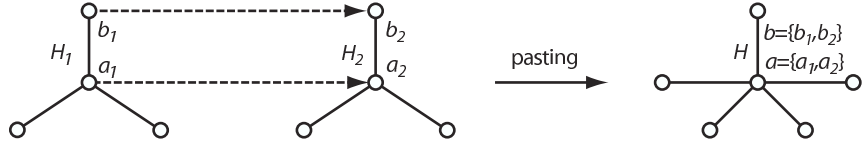}}
\caption {Pasting of two trees.} \label{example6} 
\end{figure}
}

{\begin{figure}[h!]
\centerline{\includegraphics{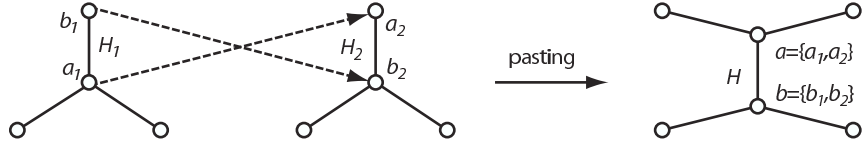}}
\caption {Another pasting of the same trees.} \label{example7} 
\end{figure}
}

In general, pasting of two partial cubes $G_1$ and $G_2$ along two isomorphic subgraphs $H_1$ and $H_2$ does not produce a partial cube even under strong assumptions about these subgraphs as the next example illustrates.

{\begin{figure}[h!]
\centerline{\includegraphics{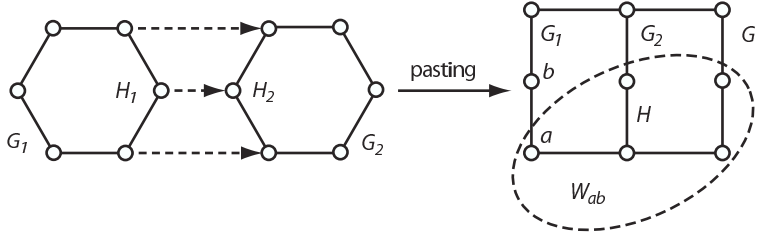}}
\caption {Pasting partial cubes $G_1$ and $G_2$.} \label{example13}
\end{figure}
}

\begin{example}
{\rm Pasting of two partial cubes $G_1=C_6$ and $G_2=C_6$ along subgraphs $H_1$ and $H_2$ is shown in Figure~\ref{example13}. The resulting graph $G$ is not a partial cube. Indeed, the semicube $W_{ab}$ is not a convex set. Note that subgraphs $H_1$ and $H_2$ are convex subgraphs of the respective partial cubes.
}
\end{example}

In this section we study two simple pastings of connected graphs together, the vertex-pasting and the edge-pasting, and show that these pastings produce partial cubes from partial cubes. We also compute the isometric and lattice dimensions of the resulting graphs.

\vtl
Let $G_1=(V_1,E_1)$ and $G_2=(V_2,E_2)$ be two connected graphs, $a_1\in V_1$, $a_2\in V_2$, and $H_1=(\{a_1\},\es),\;H_2=(\{a_2\},\es)$. Let $G$ be the graph obtained by pasting $G_1$ and $G_2$ along subgraphs $H_1$ and $H_2$. In this case we say that the graph $G$ is obtained from graphs $G_1$ and $G_2$ by {\sl vertex-pasting}. We also say that $G$ is obtained from $G_1$ and $G_2$ by {\sl identifying} vertices $a_1$ and $a_2$. Figure~\ref{vertex pasting} illustrates this construction. Note that the vertex $a=\{a_1,a_2\}$ is a cut vertex of $G$, since $G_1\cup G_2=G$ and $G_1\cap G_2=\{a\}$. (We follow our convention and identify graphs $G_1$ and $G_2$ with subgraphs of $G$.)

{\begin{figure}[h!]
\centerline{\includegraphics{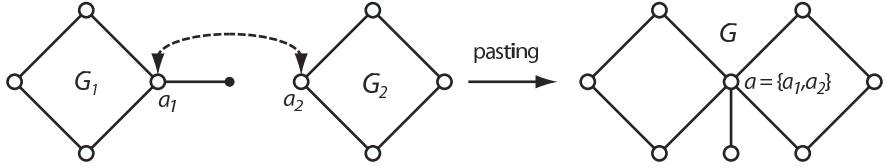}}
\caption {An example of vertex-pasting.} \label{vertex pasting} 
\end{figure}
}

In what follows we use superscripts to distinguish subgraphs of the graphs $G_1$ and $G_2$. For instance, $W^{(2)}_{ab}$ stands for the semicube of $G_2$ defined by two adjacent vertices $a,b\in V_2$.

\begin{theorem}
A graph $G=(V,E)$ obtained by vertex-pasting from partial cubes $G_1=(V_1,E_1)$ and $G_2=(V_2,E_2)$ is a partial cube.
\end{theorem}

\begin{proof}
We denote $a=\{a_1,a_2\}$ the vertex of $G$ obtained by identifying vertices $a_1\in V_1$ and $a_2\in V_2$. Clearly, $G$ is a bipartite graph. Let $xy$ be an edge of $G$. Without loss of generality we may assume that $xy\in E_1$ and $a\in W_{xy}$. Note that any path between vertices in $V_1$ and $V_2$ must go through $a$. Since $a\in W_{xy}$, we have, for any $v\in V_2$,
$$
d(v,x)=d(v,a)+d(a,x)<d(v,a)+d(a,y)=d(v,y),
$$
which implies $V_2\SB W_{xy}$ and $W_{yx}\SB V_1$. It follows that $W_{xy}=W^{(1)}_{xy}\cup V_2$ and $W_{yx}=W^{(1)}_{yx}$. The sets $W^{(1)}_{xy},\;W^{(1)}_{yx}$ and $V_2$ are convex subsets of $V$. Since $W^{(1)}_{xy}\cap V_2=\{a\}$, the set $W_{xy}=W^{(1)}_{xy}\cup V_2$ is also convex. By Theorem~\ref{DWT}(ii), the graph $G$ is a partial cube.
\end{proof}

\vtl
The vertex-pasting construction introduced above can be generalized as follows. Let $\GGG=\{G_i=(V_i,E_i)\}_{i\in J}$ be a family of connected graphs and $\AAA=\{a_i\in G_i\}_{i\in J}$ be a family of distinguished vertices of these graphs. Let $G$ be the graph obtained from the graphs $G_i$ by identifying vertices in the set $\AAA$. We say that $G$ is obtained by {\sl vertex-pasting together the graphs $G_i$} (along the set $\AAA$).

\begin{example}
{\rm Let $J=\{1,\ldots,n\}$ with $n\geq 2$, 
$$
\GGG=\{G_i=(\{a_i,b_i\},\{a_ib_i\})\}_{i\in J},\quad\text{and}\quad\AAA=\{a_i\}_{i\in J}.
$$
Clearly, each $G_i$ is $K_2$. By vertex-pasting these graphs along $\AAA$, we obtain the $n$-star graph $K_{1,n}$.
}
\end{example}

Since the star $K_{1,n}$ is a tree it can be also obtained from $K_1$ by successive vertex-pasting as in Example~\ref{tree example}.

\begin{example} \label{tree example}
{\rm Let $G_1$ be a tree and $G_2=K_2$. By vertex-pasting these graphs we obtain a new tree. Conversely, let $G$ be a tree and $v$ be its leaf. Let $G_1$ be a tree obtained from $G$ by deleting the leaf $v$. Clearly, $G$ can be obtained by vertex-pasting $G_1$ and $K_2$. It follows that any tree can obtained from the graph $K_1$ by successive vertex-pasting of copies of $K_2$ (cf. Theorem~2.3(e) in~\cite{aF95}).
}
\end{example}

Any connected graph $G$ can be constructed by successive vertex-pasting of its blocks using its {\sl block cut-vertex tree}~\cite{jB95} structure. Let $G_1$ be an endblock of $G$ with a cut vertex $v$ and $G_2$ be the union of the remaining blocks of $G$. Then $G$ can be obtained from $G_1$ and $G_2$ by vertex-pasting along the vertex $v$. It follows that any connected graph can be obtained from its blocks by successive vertex-pastings.

\vtl
Let $G=(V,E)$ be a partial cube. We recall that the isometric dimension $\dim_I(G)$ of $G$ is the cardinality of the quotient set $E\slash\hh$, where $\hh$ is Djokovi\'{c}'s equivalence relation on the set $E$ (cf. formula~(\ref{isometric dimension})).

\begin{theorem} \label{vertex isometric dim}
Let $G=(V,E)$ be a partial cube obtained by vertex-pasting together partial cubes $G_1=(V_1,E_1)$ and $G_2=(V_2,E_2)$. Then
$$
\dim_I(G)=\dim_I(G_1)+\dim_I(G_2).
$$
\end{theorem}

\begin{proof}
It suffices to prove that there are no edges $xy\in E_1$ and $uv\in E_2$ which are in Djokovi\'{c}'s relation $\hh$ with each other. Suppose that $G_1$ and $G_2$ are vertex-pasted along vertices $a_1\in E_1$ and $a_2\in E_2$ and let $a=\{a_1,a_2\}\in E$. Let $xy\in E_1$ and $uv\in E_2$ be two edges in $E$. We may assume that $u\in W_{xy}$. Since $a$ is a cut-vertex of $G$ and $u\in W_{xy}$, we have
$$
d(u,a)+d(a,x)=d(u,x)<d(u,y)=d(u,a)+d(a,y).
$$
Hence, $d(a,x)<d(a,y)$, which implies 
$$
d(v,x)=d(v,a)+d(a,x)<d(v,a)+d(a,y)=d(v,y).
$$
It follows that $v\in W_{xy}$. Therefore the edge $xy$ does not stand in the relation $\hh$ to the vertex $uv$.
\end{proof}

The next result follows immediately from the previous theorem. Note that blocks of a partial cube are partial cubes themselves.

\begin{corollary}
Let $G$ be a partial cube and $\{G_1,\ldots,G_n\}$ be the family of its blocks. Then
$$
\dim_I(G)=\sum_{i=1}^n\dim_I(G_i).
$$
\end{corollary}

\vtl
In the case of the lattice dimension of a partial cube we can claim only much weaker result than one stated in Theorem~\ref{vertex isometric dim} for the isometric dimension. We omit the proof.

\begin{theorem} \label{vertex lattice dim}
Let $G$ be a partial cube obtained by vertex-pasting together partial cubes $G_1$ and $G_2$. Then
$$
\max\{\dim_Z(G_1),\dim_Z(G_2)\}\leq\dim_Z(G)\leq\dim_Z(G_1)+\dim_Z(G_2).
$$
\end{theorem}

The following example illustrate possible cases for inequalities in Theorem~\ref{vertex lattice dim}. Let us recall that the lattice dimension of a tree with $m$ leaves is $\lceil m\slash 2\rceil$ (cf.~\cite{fH78}).

\begin{example}
{\rm The star $K_{1,6}$ can be obtained from the stars $K_{1,2}$ and $K_{1,4}$ by vertex-pasting these two stars along their centers. Clearly,
$$
\max\{\dim_Z(K_{1,2}),\dim_Z(K_{1,4})\}<\dim_Z(K_{1,6})=\dim_Z(K_{1,2})+\dim_Z(K_{1,4}).
$$

The same star $K_{1,6}$ is obtained from two copies of the star $K_{1,3}$ by vertex-pasting along their centers. We have $\dim_Z(K_{1,3})=2,\;\dim_Z(K_{1,6})=3$, so
$$
\max\{\dim_Z(K_{1,3}),\dim_Z(K_{1,3})\}<\dim_Z(K_{1,6})<\dim_Z(K_{1,3})+\dim_Z(K_{1,3}).
$$

Let us vertex-paste two stars $K_{1,3}$ along their two leaves. The resulting graph $T$ is a tree with four vertices. Therefore,
$$
\max\{\dim_Z(K_{1,3}),\dim_Z(K_{1,3})\}=\dim_Z(T)<\dim_Z(K_{1,3})+\dim_Z(K_{1,3}).
$$
}
\end{example}

\vtl
We now consider another simple way of pasting two graphs together.

\vtl
Let $G_1=(V_1,E_1)$ and $G_2=(V_2,E_2)$ be two connected graphs, $a_1b_1\in E_1$, $a_2b_2\in E_2$, and $H_1=(\{a_1,b_1\},\{a_1b_1\}),\;H_2=(\{a_2,b_2\},\{a_2b_2\})$. Let $G$ be the graph obtained by pasting $G_1$ and $G_2$ along subgraphs $H_1$ and $H_2$. In this case we say that the graph $G$ is obtained from graphs $G_1$ and $G_2$ by {\sl edge-pasting}. Figures~\ref{example6},~\ref{example7}, and~\ref{edge pasting} illustrate this construction.

{\begin{figure}[h!]
\centerline{\includegraphics{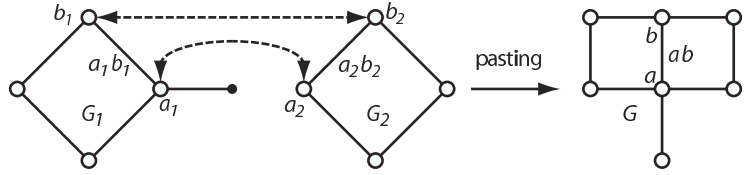}}
\caption {An example of edge-pasting.} \label{edge pasting} 
\end{figure}
}

As before, we identify the graphs $G_1$ and $G_2$ with subgraphs of the graph $G$ and denote $a=\{a_1,a_2\},\;b=\{b_1,b_2\}$ the two vertices obtained by pasting together vertices $a_1$ and $a_2$ and, respectively, $b_1$ and $b_2$. The edge $ab\in E$ is obtained by pasting together edges $a_1b_1\in E_1$ and $a_2b_2\in E_2$ (cf. Figure~\ref{edge pasting}). Then $G=G_1\cup G_2,\;V_1\cap V_2=\{a,b\}$ and $E_1\cap E_2=\{ab\}$. We use these notations in the rest of this section.

\begin{proposition} \label{edge-pasting bipartite}
A graph $G$ obtained by edge-pasting together bipartite graphs $G_1$ and $G_2$ is bipartite.
\end{proposition}

\begin{proof}
Let $C$ be a cycle in $G$. If $C\SB G_1$ or $C\SB G_2$, then the length of $C$ is even, since the graphs $G_1$ and $G_2$ are bipartite. Otherwise, the vertices $a$ and $b$ separate $C$ into two paths each of odd length. Therefore $C$ is a cycle of even length. The result follows.
\end{proof}

\vtl
The following lemma is instrumental; it describes the semicubes of the graph $G$ in terms of semicubes of graphs $G_1$ and $G_2$.

\begin{lemma} \label{main lemma}
Let $uv$ be an edge of $G$. Then
\roster
	\item[{\rm(i)}] For $uv\in E_1,\quad$$a,b\in W_{uv}\quad\imp\quad W_{uv}=W_{uv}^{(1)}\cup V_2,\;W_{vu}=W_{vu}^{(1)}$.
	\item[{\rm(ii)}] For $uv\in E_2,\quad$$a,b\in W_{uv}\quad\imp\quad W_{uv}=W_{uv}^{(2)}\cup V_1,\;W_{vu}=W_{vu}^{(2)}$.
	\item[{\rm(iii)}] $a\in W_{uv}$, $b\in W_{vu}\quad\imp\quad W_{uv}=W_{ab}$.
\endroster
\end{lemma}

{\begin{figure}[h!]
\vspace{0.1in}
\centerline{\includegraphics{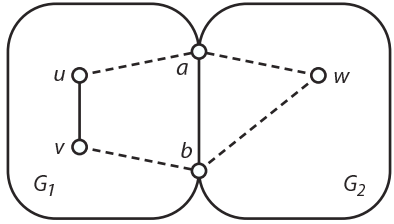}}
\caption {Edge-pasting of graphs $G_1$ and $G_2$.} \label{example9} 
\end{figure}
}

\begin{proof}
We prove parts (i) and (iii) (see Figure~\ref{example9}).

\vtl
(i) Since any path from $w\in V_2$ to $u$ or $v$ contains $a$ or $b$ and $a,b\in W_{uv}$, we have $w\in W_{uv}$. Hence, $W_{uv}=W_{uv}^{(1)}\cup V_2$ and $W_{vu}=W_{vu}^{(1)}$.

\vtl
(iii) Since $ab\,\hh\,uv$ in $G_1$, we have $W^{(1)}_{uv}=W^{(1)}_{ab}$, by Theorem~\ref{DWT}(iv). Let $w$ be a vertex in $W^{(2)}_{uv}$. Then, by the triangle inequality,
$$
d(w,u)<d(w,v)\leq d(w,b)+d(b,v)<d(w,b)+d(b,u).
$$
Since any shortest path from $w$ to $u$ contains $a$ or $b$, we have
$$
d(w,a)+d(a,u)=d(w,u).
$$
Therefore,
$$
d(w,a)+d(a,u)<d(w,b)+d(b,u).
$$
Since $ab\,\hh\,uv$ in $G_1$, we have $d(a,u)=d(b,v)$, by Theorem~\ref{F-characterization}. It follows that $d(w,a)<d(w,b)$, that is, $w\in W^{(2)}_{ab}$. We proved that $W^{(2)}_{uv}\SB W^{(2)}_{ab}$. By symmetry, $W^{(2)}_{vu}\SB W^{(2)}_{ba}$. Since two opposite semicubes form a partition of $V_2$, we have $W^{(2)}_{uv}=W^{(2)}_{ab}$. The result follows.
\end{proof}

\begin{theorem}
A graph $G$ obtained by edge-pasting together partial cubes $G_1$ and $G_2$ is a partial cube.
\end{theorem}

\begin{proof}
By Theorem~\ref{DWT}(ii) and Proposition~\ref{edge-pasting bipartite}, we need to show that for any edge $uv$ of $G$ the semicube $W_{uv}$ is a convex subset of $V$. There are two possible cases.

\vtl
(i) $uv=ab$. The semicube $W_{ab}$ is the union of semicubes $W^{(1)}_{ab}$ and $W^{(2)}_{ab}$ which are convex subsets of $V_1$ and $V_2$, respectively. It is clear that any shortest path connecting a vertex in $W^{(1)}_{ab}$ with a vertex in $W^{(2)}_{ab}$ contains vertex $a$ and therefore is contained in $W_{ab}$. Hence, $W_{ab}$ is a convex set. A similar argument proves that the set $W_{ba}$ is convex.

\vtl
(ii) $uv\not=ab$. We may assume that $uv\in E_1$. To prove that the semicube $W_{uv}$ is a convex set, we consider two cases.

\vtl
(a) $a,b\in W_{uv}$. (The case when $a,b\in W_{vu}$ is treated similarly.) By Lemma~\ref{main lemma}(i), the semicube $W_{uv}$ is the union of the semicube $W^{(1)}_{uv}$ and the set $V_2$ which are both convex sets. Any shortest path $P$ from a vertex in $V_2$ to a vertex in $W^{(1)}_{uv}$ contains either $a$ or $b$. It follows that $P\SB W_{uv}^{(1)}\cup V_2=W_{uv}$. Therefore the semicube $W_{uv}$ is convex.

\vtl
(b) $a\in W_{uv}$, $b\in W_{vu}$. (The case when $b\in W_{uv}$, $a\in W_{vu}$ is treated similarly.) By Lemma~\ref{main lemma}(ii), $W_{uv}=W_{ab}$. The result follows from part (i) of the proof.
\end{proof}

\vtl
\begin{theorem} \label{edge isometric dim}
Let $G$ be a graph obtained by edge-pasting together finite partial cubes $G_1$ and $G_2$. Then
$$
\dim_I(G)=\dim_I(G_1)+\dim_I(G_2)-1.
$$
\end{theorem}

\begin{proof}
Let $\hh$, $\hh_1$, and $\hh_2$ be Djokovi\'{c}'s relations on $E$, $E_1$, and $E_2$, respectively. By Lemma~\ref{main lemma}, for $uv,xy\in E_1$ (resp. $uv,xy\in E_2$) we have
$$
uv\,\hh\,xy\quad\Leftrightarrow\quad uv\,\hh_1 xy\quad\text{(resp. $uv\,\hh\,xy\quad\Leftrightarrow\quad uv\,\hh_2 xy$).}
$$
Let $uv\in E_1$, $xy\in E_2$, and $uv\,\hh \,xy$. Suppose that $(uv,ab)\notin\hh$. We may assume that $a,b\in W_{uv}$. By Lemma~\ref{main lemma}(i), $V_2\subset W_{uv}$, a contradiction, since $xy\in E_2$. Hence, $uv\,\hh\,xy\,\hh\,ab$. It follows that each equivalence class of the relation $\hh$ is either an equivalence class of $\hh_1$, an equivalence class of $\hh_2$ or the class containing the edge $ab$. Therefore
$$
|E\slash\hh|=|E_1\slash\hh_1|+|E_2\slash\hh_2|-1.
$$
The result follows, since the isometric dimension of a partial cube is equal to the cardinality of the set of equivalence classes of Djokovi\'{c}'s relation (formula~(\ref{isometric dimension})).
\end{proof}

We need some results about semicube graphs in order to prove an analog of Theorem~\ref{vertex lattice dim} for a partial cube obtained by edge-pasting of two partial cubes.

\begin{lemma} \label{replacement edge}
Let $G$ be a partial cube and $W_{pq}W_{uv},\;W_{qp}W_{xy}$ be two edges in the graph $\text{Sc}(G)$. Then $W_{xy}W_{uv}$ is an edge in $\text{Sc}(G)$.
\end{lemma}

\begin{proof}
By condition~(\ref{compatible2}), $W_{qp}\subset W_{uv}$ and $W_{yx}\subset W_{qp}$. Hence, $W_{yx}\subset W_{uv}$. By the same condition, $W_{xy}W_{uv}\in\text{Sc}(G)$.
\end{proof}

\vtl
As before, we identify partial cubes $G_1$ and $G_2$ with subgraphs of the partial cube $G$. Then $G_1\cup G_2=G$ and $G_1\cap G_2=(\{a,b\},\{ab\})=K_2$ (cf.~Figure~\ref{example9}).

\begin{lemma} \label{semicube edges}
Let $G$ be a partial cube obtained by edge-pasting together partial cubes $G_1$ and $G_2$. Let $W_{uv}^{(1)}W_{xy}^{(1)}$ {\rm(}resp. $W_{uv}^{(2)}W_{xy}^{(2)}${\rm)} be an edge in the semicube $\text{Sc}(G_1)$ {\rm(}resp. $\text{Sc}(G_2)${\rm)}. Then $W_{uv}W_{xy}$ is an edge in $\text{Sc}(G)$.
\end{lemma}

{\begin{figure}[h]
\centerline{\includegraphics{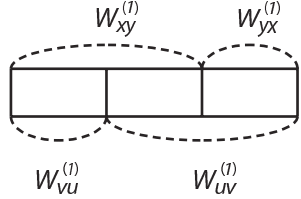}}
\caption {Semicubes forming an edge in $\text{Sc}(G_1)$.} \label{example11} 
\end{figure}
}

\begin{proof}
It suffices to consider the case of $\text{Sc}(G_1)$ (see Figure~\ref{example11}). By condition~(\ref{compatible2}), $W_{vu}^{(1)}\subset W_{xy}^{(1)}$ and $W_{yx}^{(1)}\subset W_{uv}^{(1)}$. Suppose that $a\in W_{vu}^{(1)}$ and $b\in W_{yx}^{(1)}$ (the case when $b\in W_{vu}^{(1)}$ and $a\in W_{yx}^{(1)}$ is treated similarly). Then $ab\,\hh_1 xy$ and $ab\,\hh_1 uv$. By transitivity of $\hh_1$, we have $uv\,\hh_1 xy$, a contradiction, since semicubes $W_{uv}^{(1)}$ and $W_{xy}^{(1)}$ are distinct. Therefore we may assume that, say, $a,b\in W_{uv}^{(1)}$. Then, by Lemma~\ref{main lemma}, $W_{vu}=W_{vu}^{(1)}\subset V_1$. Since $W_{vu}^{(1)}\subset W_{xy}^{(1)}\SB W_{xy}$, we have $W_{vu}\subset W_{xy}$. By condition~(\ref{compatible2}), $W_{uv}W_{xy}$ is an edge in $\text{Sc}(G)$.
\end{proof}

\vtl
\begin{lemma}
Let $M_1$ and $M_2$ be matchings in graphs $\text{Sc}(G_1)$ and $\text{Sc}(G_2)$. There is a matching $M$ in $\text{Sc}(G)$ such that
$$
|M|\geq |M_1|+|M_2|-1.
$$
\end{lemma}

\begin{proof}
By Lemma~\ref{semicube edges}, $M_1$ and $M_2$ induce matchings in $\text{Sc}(G)$ which we denote by the same symbols. The intersection $M_1\cap M_2$ is either empty or a subgraph of the empty graph with vertices $W_{ab}$ and $W_{ba}$. 

If $M_1\cap M_2$ is empty, then $M=M_1\cup M_2$ is a matching in $\text{Sc}(G)$ and the result follows.

If $M_1\cap M_2$ is an empty graph with a single vertex, say, in $M_1$, we remove from $M_1$ the edge that has this vertex as its end vertex, resulting in the matching $M'_1$. Clearly, $M=M'_1\cup M_2$ is a matching in $\text{Sc}(G)$ and $|M|=|M_1|+|M_2|-1$.

Suppose now that $M_1\cap M_2$ is the empty graph with vertices $W_{ab}$ and $W_{ba}$. Let $W_{ab}W_{uv},\;W_{ba}W_{pq}$ (resp. $W_{ab}W_{xy},\;W_{ba}W_{rs}$) be edges in $M_1$ (resp. $M_2$). By Lemma~\ref{replacement edge}, $W_{xy}W_{rs}$ is an edge in $\text{Sc}(G_2)$. Let us replace edges $W_{ab}W_{xy}$ and $W_{ba}W_{rs}$ in $M_2$ by a single edge $W_{xy}W_{rs}$, resulting in the matching $M'_2$. Then $M=M_1\cup M'_2$ is a matching in $\text{Sc}(G)$ and $|M|=|M_1|+|M_2|-1$.
\end{proof}

\vtl
\begin{corollary}
Let $M_1$ and $M_2$ be maximum matchings in $\text{Sc}(G_1)$ and $\text{Sc}(G_2)$, respectively, and $M$ be a maximum matching in $\text{Sc}(G)$. Then
\begeq \label{M>=M1+M2-1}
|M|\geq |M_1|+|M_2|-1.
\edeq
\end{corollary}

\vtl
By Theorem~\ref{lattice=isometric-|M|}, we have
$$
\dim_I(G_1)=\dim_Z(G_1)+|M_1|,\quad \dim_I(G_2)=\dim_Z(G_2)+|M_2|,
$$
and
$$
\dim_I(G)=\dim_Z(G)+|M|,
$$
where $M_1$ and $M_2$ are maximum matchings in $\text{Sc}(G_1)$ and $\text{Sc}(G_2)$, respectively, and $M$ is a maximum matching in $\text{Sc}(G)$. Therefore, by Theorem~\ref{edge isometric dim} and~(\ref{M>=M1+M2-1}), we have the following result (cf. Theorem~\ref{vertex lattice dim}).

\begin{theorem} \label{edge lattice dim}
Let $G$ be a partial cube obtained by edge-pasting from partial cubes $G_1$ and $G_2$. Then
$$
\max\{\dim_Z(G_1),\dim_Z(G_2)\}\leq\dim_Z(G)\leq\dim_Z(G_1)+\dim_Z(G_2).
$$
\end{theorem}

\vtl
\begin{example}
{\rm Let us consider two edge-pastings of the stars $G_1=K_{1,3}$ and $G_2=K_{1,3}$ of lattice dimension $2$ shown in figures~\ref{example6} and~\ref{example7}. In the first case the resulting graph is the star $G=K_{1,5}$ of lattice dimension $3$.
Then we have
$$
\max\{\dim_Z(G_1),\dim_Z(G_2)\}<\dim_Z(G)<\dim_Z(G_1)+\dim_Z(G_2).
$$
In the second case the resulting graph is a tree with $4$ leaves. Therefore,
$$
\max\{\dim_Z(G_1),\dim_Z(G_2)\}=\dim_Z(G)<\dim_Z(G_1)+\dim_Z(G_2).
$$

Let $c_1a_1$ and $c_2a_2$ be edges of stars $G_1=K_{1,4}$ and $G_2=K_{1,4}$ (each of which has lattice dimension $2$), where $c_1$ and $c_2$ are centers of the respective stars. Let us edge-paste these two graphs by identifying $c_1$ with $c_2$ and $a_1$ with $a_2$, respectively. The resulting graph $G$ is the star $K_{1,7}$ of lattice dimension $4$. Thus,
$$
\max\{\dim_Z(G_1),\dim_Z(G_2)\}\leq\dim_Z(G)=\dim_Z(G_1)+\dim_Z(G_2).
$$
}
\end{example}

\section{Expansions and contractions of partial cubes} \label{S:expansions}

The graph expansion procedure was introduced by Mulder in~\cite{hM80}, where it is shown that a graph is a median graph if and only if it can be obtained from $K_1$ by a sequence of convex expansions (see also~\cite{wI00}). A similar result for partial cubes was established in~\cite{vC88} (see also~\cite{vC94}) as a corollary to a more general result concerning isometric embeddability into Hamming graphs; it was also established in~\cite{kF93} in the framework of oriented matroids theory. 

\vtl
In this section we investigate properties of (isometric) expansion and contraction operations and, in particular, prove in two different ways that a graph is a partial cube if and only if it can be obtained from the graph $K_1$ by a sequence of expansions.

\vtl
A remark about notations is in order. In the product $\{1,2\}\times (V_1\cup V_2)$,
we denote $V'_i=\{i\}\times V_i$ and $x^i=(i,x)$ for $x\in V_i$, where $i,j=1,2$.

\begin{definition} \label{D:expansion}
{\rm Let $G=(V,E)$ be a connected graph, and let $G_1=(V_1,E_1)$ and $G_2=(V_2,E_2)$ be two isometric subgraphs of $G$ such that $G=G_1\cup G_2$. The {\sl expansion of $G$ with respect to $G_1$ and $G_2$} is the graph $G'=(V',E')$ constructed as follows from $G$ (see Figure~\ref{example14}):
\roster
	\item[(i)] $V'=V_1+V_2=V'_1\cup V'_2$;
	\item[(ii)] $E'=E_1+E_2+M$, where $M$ is the matching $\bigcup_{x\in V_1\cap V_2}\{x^1x^2\}$.
\endroster
In this case, we also say that $G$ is a {\sl contraction} of $G'$.
}
\end{definition}

{\begin{figure}[h!]
\centerline{\includegraphics{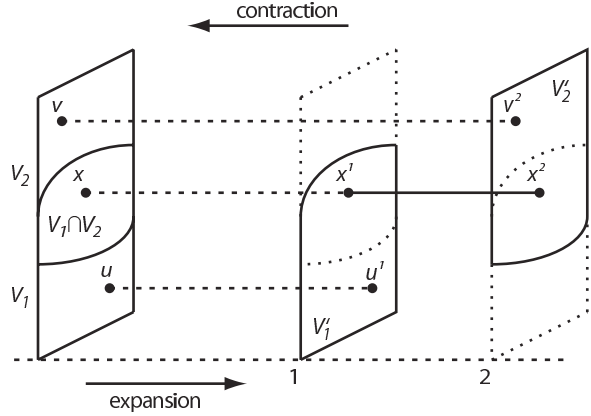}}
\caption {Expansion/contraction processes.} \label{example14} 
\end{figure}
}

It is clear that the graphs $G_1$ and $\langle V'_1\rangle$ are isomorphic, as well as the graphs $G_2$ and $\langle V'_2\rangle$.

\vtl
We define a {\sl projection} $p:V'\rightarrow V$ by $p(x^i)=x$ for $x\in V$. Clearly, the restriction of $p$ to $V'_1$ is a bijection $p_1:V'_1\rightarrow V_1$ and its restriction to $V'_2$ is a bijection $p_2:V'_2\rightarrow V_2$. These bijections define isomorphisms $\langle V'_1\rangle\rightarrow G_1$ and $\langle V'_2\rangle\rightarrow G_2$.

\vtl
Let $P'$ be a path in $G'$. The vertices of $G$ obtained from the vertices in $P'$ under the projection $p$ define a walk $P$ in $G$; we call this walk $P$ the {\sl projection} of the path $P'$.
It is clear that
\begeq \label{l(P)=l(P')}
\ell(P)=\ell(P'),\quad\text{if $P'\SB\langle V'_1\rangle$ or $P'\SB\langle V'_2\rangle$.}
\edeq
In this case, $P$ is a path in $G$ and either $P=p_1(P')$ or $P=p_2(P')$. On the other hand,
\begeq \label{l(P)<l(P')}
\ell(P)<\ell(P'),\quad\text{if $P'\cap\langle V'_1\rangle\not=\es$ and $P'\cap\langle V'_2\rangle\not=\es$,}
\edeq
and $P$ is not necessarily a path.

\vtl
We will frequently use the results of the following lemma in this section.

\begin{lemma} \label{d=d+1}
{\rm(i)} For $u^1,v^1\in V'_1$, any shortest path $P_{u^1v^1}$ in $G'$ belongs to $\langle V'_1\rangle$ and its projection $P_{uv}=p_1(P_{u^1v^1})$ is a shortest path in $G$. Accordingly,
$$
d_{G'}(u^1,v^1)=d_G(u,v)
$$
and $\langle V'_1\rangle$ is a convex subgraph of $G'$. A similar statement holds for $u^2,v^2\in V'_2$.

\vtl
{\rm(ii)} For $u^1\in V'_1$ and $v^2\in V'_2$,
$$
d_{G'}(u^1,v^2)=d_G(u,v)+1.
$$
Let $P_{u^1v^2}$ be a shortest path in $G'$. There is a unique edge $x^1x^2\in M$ such that $x^1,x^2\in P_{u^1v^2}$ and the sections $P_{u^1x^1}$ and $P_{x^2v^2}$ of the path $P_{u^1v^2}$ are shortest paths in $\langle V'_1\rangle$ and $\langle V'_2\rangle$, respectively. The projection $P_{uv}$ of $P_{u^1v^2}$ in $G'$ is a shortest path in $G$.
\end{lemma}

\begin{proof}
(i) Let $P_{u^1v^1}$ be a path in $G'$ that intersects $V'_2$. Since  $\langle V_1\rangle$ is an isometric subgraph of $G$, there is a path $P_{uv}$ in $G$ that belongs to $\langle V_1\rangle$. Then $p_1^{-1}(P_{uv})$ is a path in $\langle V'_1\rangle$ of the same length as $P_{uv}$. By~(\ref{l(P)=l(P')}) and~(\ref{l(P)<l(P')}), 
$$
\ell(p_1^{-1}(P_{uv}))<\ell(P_{u^1v^1}).
$$
Therefore any shortest path $P_{u^1v^1}$ in $G'$ belongs to $\langle V'_1\rangle$. The result follows.

\vtl
(ii) Let $P_{u^1v^2}$ be a shortest path in $G'$ and $P_{uv}$ be its projection to $V$. By~(\ref{l(P)<l(P')}),
$$
d_{G'}(u^1,v^2)=\ell(P_{u^1v^2})>\ell(P_{uv})\geq d_G(u,v).
$$
Since there is no edge of $G$ joining vertices in $V_1\setminus V_2$ and $V_2\setminus V_1$, a shortest path in $G$ from $u$ to $v$ must contain a vertex $x\in V_1\cap V_2$. Since $G_1$ and $G_2$ are isometric subgraphs, there are shortest paths $P_{ux}$ in $G_1$ and $P_{xv}$ in $G_2$ such that their union is a shortest path from $u$ to $v$. Then, by the triangle inequality and part (i) of the proof, we have (cf. Figure~\ref{example14})
$$
d_{G'}(u^1,v^2)\leq d_{G'}(u^1,x^1)+d_{G'}(x^1,x^2)+d_{G'}(x^2,v^2)=d_G(u,v)+1.
$$
The last two displayed formulas imply $d_{G'}(u^1,v^2)=d_G(u,v)+1$.

Since $u^1\in V'_1$ and $v^2\in V'_2$ the path $P_{u^1v^2}$ must contain an edge, say $x^1x^2$, in $M$. Since this path is a shortest path in $G'$, this edge is unique. Then the sections $P_{u^1x^1}$ and $P_{x^2v^2}$ of $P_{u^1v^2}$ are shortest paths in $\langle V'_1\rangle$ and $\langle V'_2\rangle$, respectively. Clearly, $P_{uv}$ is a shortest path in $G$.
\end{proof}

\vtl
Let $a^1a^2$ be an edge in the matching $M=\cup_{x\in V_1\cap V_2}\{x^1x^2\}$. This edge defines five fundamental sets (cf. Section~\ref{S:fundamental sets}): the semicubes $W_{a^1a^2}$ and $W_{a^2a^1}$, the sets of vertices $U_{a^1a^2}$ and $U_{a^2a^1}$, and the set of edges $F_{a^1a^2}$. The next theorem follows immediately from Lemma~\ref{d=d+1}. It gives a hint to a connection between the expansion process and partial cubes.

\begin{theorem} \label{semicubes in expansion}
Let $G'$ be an expansion of a connected graph $G$ and notations are chosen as above. Then
\roster
	\item[{\rm(i)}] $W_{a^1a^2}=V'_1$ and $W_{a^2a^1}=V'_2$ are convex semicubes of $G'$.
	\item[{\rm(ii)}] $F_{a^1a^2}=M$ defines an isomorphism between induced subgraphs $\langle U_{a^1a^2}\rangle$ and $\langle U_{a^2a^1}\rangle$, which are isomorphic to the subgraph $G_1\cap G_2$.
\endroster
\end{theorem}

The result of Theorem~\ref{semicubes in expansion} justifies the following constructive definition of the contraction process.

\begin{definition} \label{D:contraction}
{\rm Let $ab$ be an edge of a connected graph $G'=(V',E')$ such that
\roster
	\item[(i)] semicubes $W_{ab}$ and $W_{ba}$ are convex and form a partition of $V'$;
	\item[(ii)] the set $F_{ab}$ is a matching and defines an isomorphism between subgraphs $\langle U_{ab}\rangle$ and $\langle U_{ba}\rangle$.
\endroster
A graph $G$ obtained from the graphs $\langle W_{ab}\rangle$ and $\langle W_{ba}\rangle$ by pasting them along subgraphs $\langle U_{ab}\rangle$ and $\langle U_{ba}\rangle$ is said to be a {\sl contraction} of the graph $G'$.
}
\end{definition}

\begin{remark}
{\rm If $G'$ is bipartite, then semicubes $W_{ab}$ and $W_{ba}$ form a partition of its vertex set. Then, by Theorem~\ref{theorem matching=isomorphism}, condition (i) implies condition (ii). Thus any pair of opposite convex semicubes in a connected bipartite graph defines a contraction of this graph.
}
\end{remark}

By Theorem~\ref{semicubes in expansion}, a graph is a contraction of its expansion. It is not difficult to see that any connected graph is also an expansion of its contraction.

\vtl
The following three examples give geometric illustrations for the expansion and contraction procedures.

\begin{example} \label{isometric not convex}
{\rm Let $a$ and $b$ be two opposite vertices in the graph $G=C_4$. Clearly, the two distinct paths $P_1$ and $P_2$ from $a$ to $b$ are isometric subgraphs of $G$ defining an expansion $G'=C_6$ of $G$ (see Figure~\ref{example15}). Note that $P_1$ and $P_2$ are not convex subsets of $V$.
}
\end{example}

{\begin{figure}[h!]
\centerline{\includegraphics[scale=0.9]{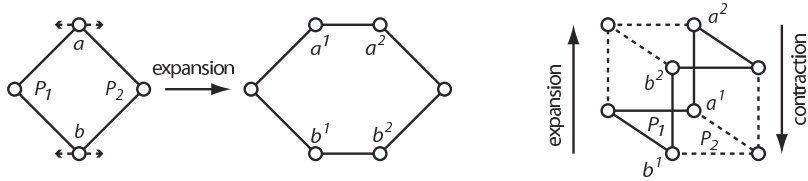}}
\caption {An expansion of the cycle $C_4$.} \label{example15} 
\end{figure}
}

{\begin{figure}[h!]
\centerline{\includegraphics[scale=0.9]{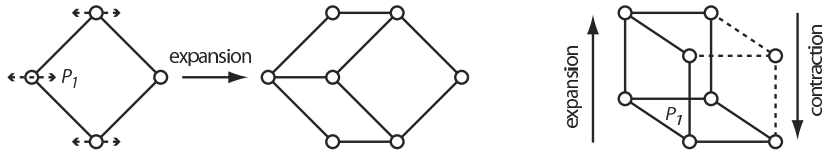}}
\caption {Another isometric expansion of the cycle $C_4$.} \label{example16} 
\end{figure}
}

\begin{example}
{\rm Another isometric expansion of the graph $G=C_4$ is shown in Figure~\ref{example16}. Here, the path $P_1$ is the same as in the previous example and $G_2=G$.
}
\end{example}

\begin{example}
{\rm Lemma~\ref{d=d+1} claims, in particular, that the projection of a shortest path in an extension $G'$ of a graph $G$ is a shortest path in $G$. Generally speaking, the converse is not true. Consider the graph $G$ shown in Figure~\ref{example17} and two paths in $G$:
$$
V_1=abcef\quad\text{and}\quad V_2=bde.
$$
The graph $G'$ in Figure~\ref{example17} is the convex expansion of $G$ with respect to $V_1$ and $V_2$. The path $abde\!f$ is a shortest path in $G$; it is not a projection of a shortest path in $G'$.
}
\end{example}

{\begin{figure}[h!]
\centerline{\includegraphics[scale=0.9]{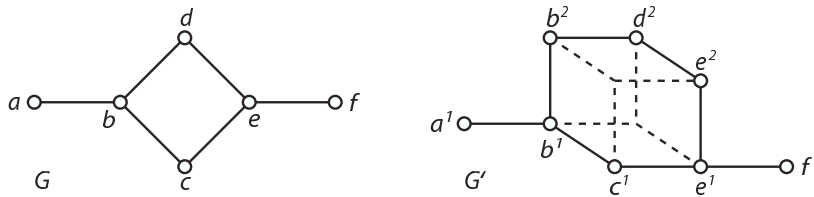}}
\caption {A shortest path which is not a projection of a shortest path.} \label{example17} 
\end{figure}
}

One can say that, in the case of finite partial cubes, the contraction procedure is defined by an orthogonal projection of a hypercube onto one of its facets.

\vtl
By Theorem~\ref{semicubes in expansion}, the sets $V'_1$ and $V'_2$ are opposite semicubes of the graph $G'$ defined by edges in $M$. Their projections are the sets $V_1$ and $V_2$ which are not necessarily semicubes of $G$. For other semicubes in $G'$ we have the following result.

\begin{lemma} \label{projection of semicube}
For any two adjacent vertices $u,v\in V$,
$$
W_{u^iv^i}=p^{-1}(W_{uv})\quad\text{for $u,v\in V_i$ and $i=1,2$.}
$$
\end{lemma}

\begin{proof}
By Lemma~\ref{d=d+1},
$$
d_{G'}(x^j,u^i)<d_{G'}(x^j,v^i)\quad\eq\quad d_G(x,u)<d_G(x,v)
$$
for $x\in V$ and $i,j=1,2$. The result follows.
\end{proof}

\vtl
\begin{corollary}
If $uv$ is an edge of $G_1\cap G_2$, then $W_{u^1v^1}=W_{u^2v^2}$.
\end{corollary}

The following lemma is an immediate consequence of Lemma~\ref{d=d+1}. We shall use it implicitly in our arguments later.

\begin{lemma}
Let $u,v\in V_1$ and $x\in V_1\cap V_2$. Then
$$
x^1\in W_{u^1v^1}\quad\Leftrightarrow\quad x^2\in W_{u^1v^1}.
$$
The same result holds for semicubes in the form $W_{u^2v^2}$.
\end{lemma}

Generally speaking, the projection of a convex subgraph of $G'$ is not a convex subgraph of $G$. For instance, the projection of the convex path $b^2d^2e^2$ in Figure~\ref{example17} is the path $bde$ which is not a convex subgraph of $G$. On the other hand, we have the following result.

\begin{theorem} \label{convex projection}
Let $G'=(V',E')$ be an expansion of a graph $G=(V,E)$ with respect to subgraphs $G_1=(V_1,E_1)$ and $G_2=(V_2,E_2)$. The projection of a convex semicube of $G'$ different from $\langle V'_1\rangle$ and $\langle V'_2\rangle$ is a convex semicube of $G$.
\end{theorem}

\begin{proof}
It suffices to consider the case when $W_{uv}=p(W_{u^1v^1})$ for $u,v\in V_1$ (cf. Theorem~\ref{projection of semicube}). Let $x,y\in W_{uv}$ and $z\in V$ be a vertex such that
$$
d_G(x,z)+d_G(z,y)=d_G(x,y).
$$
We need to show that $z\in W_{uv}$.

{\begin{figure}[h!]
\centerline{\includegraphics{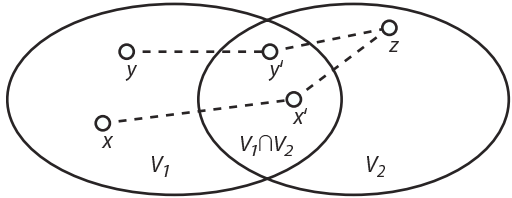}}
\caption {A shortest path from $x$ to $y$.} \label{example18} 
\end{figure}
}

(i) $x,y\in V_1$ (the case when $x,y\in V_2$ is treated similarly). Suppose that $z\in V_1$. Then $x^1,y^1,z^1\in V'_1$ and, by Lemma~\ref{d=d+1},
$$
d_{G'}(x^1,z^1)+d_{G'}(z^1,y^1)=d_{G'}(z^1,y^1).
$$
Since $x^1,y^1\in W_{u^1v^1}$ and $W_{u^1v^1}$ is convex, $z^1\in W_{u^1v^1}$. Hence, $z\in W_{uv}$.

Suppose now that $z\in V_2\setminus V_1$. Consider a shortest path $P_{xy}$ in $G$ from $x$ to $y$ containing $z$. This path contains vertices $x',y'\in V_1\cap V_2$ such that (see Figure~\ref{example18})
$$
d_G(x,x')+d_G(x',z)=d_G(x,z)\quad\text{and}\quad d_G(y,y')+d_G(y',z)=d_G(y,z).
$$
Since $P_{xy}$ is a shortest path in $G$, we have
$$
d_G(x,x')+d_G(x',y)=d_G(x,y),\quad d_G(x,y')+d_G(y',y)=d_G(x,y),
$$
and
$$
d_G(x',z)+d_G(z,y')=d_G(x',y').
$$
Since $x,x',y\in V_1$, we have $x^1,x'^1,y^1\in V'_1$. Because $x^1,y^1\in W_{u^1v^1}$ and $W_{u^1v^1}$ is convex, $x'^1\in W_{u^1v^1}$. Hence, $x'\in W_{uv}$ and, similarly, $y'\in W_{uv}$. Since $x'^2,y'^2,z^2\in V'_2$ and $W_{u^1v^1}$ is convex, $z^2\in W_{u^1v^1}$. Hence, $z\in W_{uv}$.

\vtl
(ii) $x\in V_1\setminus V_2$ and $y\in V_2\setminus V_1$. We may assume that $z\in V_1$. By Lemma~\ref{d=d+1},
\begin{align*}
d_{G'}(x^1,y^2)&=d_G(x,y)+1=d_G(x,z)+d_G(z,y)+1\\
&=d_{G'}(x^1,z^1)+d_{G'}(z^1,y^2).
\end{align*}
Since $x^1,y^2\in W_{u^1v^1}$ and $W_{u^1v^1}$ is convex, $z^1\in W_{u^1v^1}$. Hence, $z\in W_{uv}$.
\end{proof}

By using the results of Lemma~\ref{d=d+1}, it is not difficult to show that the class of connected bipartite graphs is closed under the expansion and contraction operations. The next theorem establishes this result for the class of partial cubes.

\begin{theorem} \label{expansion/contraction of partial cube}
{\rm(i)} 
An expansion $G'$ of a partial cube $G$ is a partial cube.
\vtl
{\rm(ii)} A contraction $G$ of a partial cube $G'$ is a partial cube.
\end{theorem}

\begin{proof}
(i) Let $G=(V,E)$ be a partial cube and $G'=(V',E')$ be its expansion with respect to isometric subgraphs $G_1=(V_1,E_1)$ and $G_2=(V_2,E_2)$. By Theorem~\ref{DWT}(ii), it suffices to show that the semicubes of $G'$ are convex.

\vtl
By Lemma~\ref{semicubes in expansion}, the semicubes $\langle V'_1\rangle$ and $\langle V'_2\rangle$ are convex, so we consider a semicube in the form $W_{u^1v^1}$ where $uv\in E_1$ (the other case is treated similarly). Let $P_{x'y'}$ be a shortest path connecting two vertices in $W_{u^1v^1}$ and $P_{xy}$ be its projection to $G$. By Theorem~\ref{projection of semicube}, $x,y\in W_{uv}$ and, by Lemma~\ref{d=d+1}, $P_{xy}$ is a shortest path in $G$. Since $W_{uv}$ is convex, $P_{xy}$ belongs to $W_{uv}$. Let $z'$ be a vertex in $P_{x'y'}$ and $z=p(z')\in P_{xy}$. By Lemma~\ref{d=d+1},
$$
d_G(z,u)<d_G(z,v)\quad\Rightarrow\quad d_{G'}(z',u^1)\leq d_{G'}(z',v^1).
$$
Since $G'$ is a bipartite graph, $d_{G'}(z',u^1)<d_{G'}(z',v^1)$. Hence, $P_{x'y'}\SB W_{u^1v^1}$, so $W_{u^1v^1}$ is convex.
\vtl
(ii) Let $G=(V,E)$ be a contraction of a partial cube $G'=(V',E')$. By Theorem~\ref{DWT}, we need to show that the semicubes of $G$ are convex. By Theorem~\ref{projection of semicube}, all semicubes of $G$ are projections of semicubes of $G'$ distinct from $\langle V'_1\rangle$ and $\langle V'_2\rangle$. By Theorem~\ref{convex projection}, the semicubes of $G$ are convex.
\end{proof}

\begin{corollary} \label{finite expansion/contraction}
{\rm(i)} A finite connected graph is a partial cube if and only if it can be obtained from $K_1$ by a sequence of expansions.

{\rm(ii)} The number of expansions needed to produce a partial cube $G$ from $K_1$ is $\dim_I(G)$.
\end{corollary}

\begin{proof}
(i) Follows immediately from Theorem~\ref{expansion/contraction of partial cube}.

\vtl
(ii) Follows from theorems~\ref{projection of semicube} and~\ref{Djokovic dimension} (see the discussion in Section~\ref{S:dimensions} just before Theorem~\ref{X-dimension} ).
\end{proof}

\vtl
The processes of expansion and contraction admit useful descriptions in the case of partial cubes on a set. Let $G=(V,E)$ be a partial cube on a set $X$, that is an isometric subgraph of the hypercube $\HHH(X)$. Then it is induced by some wg-family $\FFF$ of finite subsets of $X$ (cf. Theorem~\ref{wg=finite wg}). We may assume (see Section~\ref{S:dimensions}) that $\cap\,\FFF=\es$ and $\cup\,\FFF=X$.

\vtl
In what follows we present proofs of the results of Theorem~\ref{expansion/contraction of partial cube} and Corollary~\ref{finite expansion/contraction} given in terms of wg-families of sets.

\vtl
The expansion process for a partial cube $G$ on $X$ can be described as follows: Let $\FFF_1$ and $\FFF_2$ be wg-families of finite subsets of $X$ such that $\FFF_1\cap\FFF_2\not=\es$, $\FFF_1\cup\FFF_2=\FFF$, and the distance between any two sets $P\in\FFF_1\setminus\FFF_2$ and $Q\in\FFF_2\setminus\FFF_1$ is greater than one. Note that $\langle\FFF_1\rangle$ and $\langle\FFF_2\rangle$ are partial cubes, $\langle\FFF_1\rangle\cap\langle\FFF_2\rangle\not=\es$, and $\langle\FFF_1\rangle\cup\langle\FFF_2\rangle=\langle\FFF\rangle=G$. Let $X'=X+\{p\}$, where $p\notin X$, and
$$
\FFF'_2=\{Q+\{p\}: Q\in\FFF_2\},\quad\FFF'=\FFF_1\cup\FFF'_2.
$$
It is quite clear that the graphs $\langle\FFF'_2\rangle$ and $\langle\FFF_2\rangle$ are isomorphic and the graph $G'=\langle\FFF'\rangle$ is an isometric expansion of the graph $G$.

\begin{theorem}
An expansion of a partial cube is a partial cube.
\end{theorem}

\begin{proof}
We need to verify that $\FFF'$ is a wg-family of finite subsets of $X'$. By Theorem~\ref{local wg}, it suffices to show that the distance between any two adjacent sets in $\FFF'$ is $1$. It is obvious if each of these two sets belong to one of the families $\FFF_1$ or $\FFF'_2$. Suppose that $P\in\FFF_1$ and $Q+\{p\}\in\FFF'_2$ are adjacent, that is, for any $S\in\FFF'$ we have
\begeq \label{implication}
P\cap(Q+\{p\})\SB S\SB P\cup(Q+\{p\})\quad\Rightarrow\quad S=P\text{~or~}S=Q+\{p\}.
\edeq

If $Q\in\FFF_1$, then
$$
P\cap(Q+\{p\})\SB Q\SB P\cup(Q+\{p\}),
$$
since $p\notin P$. By~(\ref{implication}), $Q=P$ implying $d(P,Q+\{p\})=1$.

\vtl
If $Q\in\FFF_2\setminus\FFF_1$, there is $R\in\FFF_1\cap\FFF_2$ such that
$$
d(P,R)+d(R,Q)=d(P,Q),
$$
since $\FFF$ is well graded. By Theorem~\ref{lattice=metric betweenness},
$$
P\cap Q\SB R\SB P\cup Q,
$$
which implies
$$
P\cap(Q+\{p\})\SB R+\{p\}\SB P\cup(Q+\{p\}).
$$
By~(\ref{implication}), $R+\{p\}=Q+\{p\}$, a contradiction.
\end{proof}

\vtl
It is easy to recognize the fundamental sets (cf. Section~\ref{S:fundamental sets}) in an isometric expansion $G'$ of a partial cube $G=\langle\FFF\rangle$.
Let $P\in\FFF_1\cap\FFF_2$ and $Q=P+\{p\}\in\FFF'_2$ be two vertices defining an edge in $G'$ according to Definition~\ref{D:expansion}(ii). Clearly, the families $\FFF_1$ and $\FFF'_2$ are the semicubes $W_{PQ}$ and $W_{QP}$ of the graph $G'$ (cf. Lemma~\ref{points=semicubes}) and therefore are convex subsets of $\FFF'$. The set $F_{PQ}$ is the set of edges defined by $p$ as in Lemma~\ref{points=semicubes}. In addition, $U_{PQ}=\FFF_1\cap\FFF_2$ and $U_{QP}=\{R+\{p\}: R\in\FFF_1\cap\FFF_2\}$.

\vtl
Let $G$ be a partial cube induced by a wg-family $\FFF$ of finite subsets of a set $X$. As before, we assume that $\cap\,\FFF=\es$ and $\cup\,\FFF= X$. Let $PQ$ be an edge of $G$. We may assume that $Q=P+\{p\}$ for some $p\notin P$. Then (see Lemma~\ref{points=semicubes})
$$
W_{PQ}=\{R\in\FFF: p\notin R\}\quad\text{and}\quad W_{QP}=\{R\in\FFF: p\in R\}.
$$

Let $X'=X\setminus\{p\}$ and $\FFF'=\{R\setminus\{p\}: R\in\FFF\}$. It is clear that the graph $G'$ induced by the family $\FFF'$ is isomorphic to the contraction of $G$ defined by the edge $PQ$. Geometrically, the graph $G'$ is the orthogonal projection of the graph $G$ along the edge $PQ$ (cf. figures~\ref{example15} and~\ref{example16}).

\begin{theorem}
{\rm(i)} A contraction $G'$ of a partial cube $G$ is a partial cube. 

{\rm(ii)} If $G$ is finite, then $\dim_I(G')=\dim_I(G)-1$.
\end{theorem}

\begin{proof}
(i) For $p\in X$ we define $\FFF_1=\{R\in\FFF: p\notin R\}$, $\FFF_2=\{R\in\FFF: p\in R\}$, and $\FFF'_2=\{R\setminus\{p\}\in\FFF: p\in R\}$. Note that $\FFF_1$ and $\FFF_2$ are semicubes of $G$ and $\FFF'_2$ is isometric to $\FFF_2$. Hence, $\FFF_1$ and $\FFF'_2$ are wg-families of finite subsets of $X'$. We need to prove that $\FFF'=\FFF_1\cup\FFF'_2$ is a wg-family. By Theorem~\ref{local wg}, it suffices to show that $d(P,Q)=1$ for any two adjacent sets $P,Q\in\FFF'$. This is true if $P,Q\in\FFF_1$ or $P,Q\in\FFF'_2$, since these two families are well graded. For $P\in\FFF_1\setminus\FFF'_2$ and $Q\in\FFF'_2\setminus\FFF_1$, the sets $P$ and $Q+\{p\}$ are not adjacent in $\FFF$, since  $\FFF$ is well graded and $Q\notin\FFF$. Hence there is $R\in\FFF_1$ such that
$$
P\cap(Q+\{p\})\SB R\SB P\cup(Q+\{p\})
$$
and $R\not= P$. Since $p\notin R$, we have
$$
P\cap Q\SB R\SB P\cup Q.
$$
Since $R\not= P$ and $R\not= Q$, the sets $P$ and $Q$ are not adjacent in $\FFF'$. The result follows.

\vtl
(ii) If $G$ is a finite partial cube, then, by Theorem~\ref{X-dimension},
$$
\dim_I(G')=|X'|=|X|-1=\dim_I(G)-1.
$$
\end{proof}

\section{Conclusion}

The paper focuses on two themes of a rather general mathematical nature.

\vtl
1. {\sl The characterization problem.} It is a common practice in mathematics to characterize a particular class of object in different terms. We present new characterizations of the classes of bipartite graphs and partial cubes, and give new proofs for known characterization results.

\vtl
2. {\sl Constructions.} The problem of constructing new objects from old ones is a standard topic in many branches of mathematics. For the class of partial cubes, we discuss operations of forming the Cartesian product, expansion and contraction, and pasting. It is shown that the class of partial cubes is closed under these operations.

\vtl
Because partial cubes are defined as graphs isometrically embeddable into hypercubes, the theory of partial cubes has a distinctive geometric flavor. The three main structures on a graph---semicubes and Djokovi\'{c}'s and Winkler's relations---are defined in terms of the metric structure on a graph. One can say that this theory is a branch of discrete metric geometry. Not surprisingly, geometric structures play an important role in our treatment of the characterization and construction problems.


\end{document}